\documentclass[10pt]{article}
\usepackage{amssymb}
\usepackage{epsfig}
\usepackage[french,english]{babel}
\usepackage{graphicx}
\usepackage[latin1]{inputenc}
{\catcode `\@=11 \global\let\AddToReset=\@addtoreset}
\AddToReset{equation}{section}

\newtheorem{@definition}{\sc Definition}[section]

\newtheorem{@remark}{\sc Remark}[section]

\newtheorem{@example}{\sc Example}[section]

\newcommand{\beqn}{\begin{displaymath}}
\newcommand{\eeqn}{\end{displaymath}}
\newcommand{\beq}{\begin{equation}}  % numbered (single equation)
\newcommand{\eeq}{\end{equation}}

\def\mathsf{\bf}
\def\N{\mathbb{N}}

\def\L{\mathbb{L}}
\def\R{\mathbb{R}}
\def\Z{\mathbb{Z}}

\def\i{\mathrm i}
\def\d{\mathrm d}
\def\e{\mathrm e}

\def\E{\mathrm E}
\def\P{\mathrm P}

\def\text{\mbox}

\def\1{{\bf 1}}

\newcommand{\mbf}[1]{\mbox{\boldmath $#1$}}

\newcommand{\HH}{\mathfrak{H}}

\newcommand{\Cov}{\mbox{Cov}}

\def\limiteloin{\renewcommand{\arraystretch}{0.5}
\begin{array}[t]{c}
\stackrel{{\cal D}}{\longrightarrow} \\
{\scriptstyle n \rightarrow\infty}
\end{array}\renewcommand{\arraystretch}{1}}

\def\limiten{\renewcommand{\arraystretch}{0.5}
\begin{array}[t]{c}
\stackrel{}{\longrightarrow} \\
{\scriptstyle n\rightarrow\infty}
\end{array}\renewcommand{\arraystretch}{1}}

\def\limiteK{\renewcommand{\arraystretch}{0.5}
\begin{array}[t]{c}
\stackrel{}{\longrightarrow} \\
{\scriptstyle K \rightarrow\infty}
\end{array}\renewcommand{\arraystretch}{1}}

\def\limitepsn{\renewcommand{\arraystretch}{0.5}
\begin{array}[t]{c}
\stackrel{a.s.}{\longrightarrow} \\
{\scriptstyle n\rightarrow\infty}
\end{array}\renewcommand{\arraystretch}{1}}

\def\limitet{\renewcommand{\arraystretch}{0.5}
\begin{array}[t]{c}
\stackrel{}{\longrightarrow} \\
{\scriptstyle t\rightarrow\infty}
\end{array}\renewcommand{\arraystretch}{1}}

\def\limitet0{\renewcommand{\arraystretch}{0.5}
\begin{array}[t]{c}
\stackrel{}{\longrightarrow} \\
{\scriptstyle t\rightarrow 0}
\end{array}\renewcommand{\arraystretch}{1}}

\newtheorem{thm}{Theorem}
\newtheorem{rem}{Remark}
\newtheorem{cor}{Corollary}
\newtheorem{lem}{Lemma}
\newtheorem{prop}{Proposition}

\def\Cov{\mathrm{Cov}}

\begin{document}
%\begin{opening}
\title{\bf Moment bounds and central limit theorems for Gaussian
subordinated arrays}
\author{\centerline{Jean-Marc Bardet${}^1$ and Donatas~Surgailis${}^{2,}$\thanks
{Supported by a grant (No.\ MIP-11155) from the Research Council of Lithuania }
} \\
\small {\tt bardet@univ-paris1.fr}, \small {\tt
donatas.surgailis@mii.vu.lt}\\
~\\
{\small $^1$ SAMM, Universit\'e Paris 1, 90 rue de Tolbiac, 75013 Paris, FRANCE}\\
{\small $^{2}$Institute of Mathematics and Informatics,  Vilnius University, Akademijos
4, 08663 Vilnius, LITHUANIA}}
%\runningauthor{J.M. Bardet, D. Surgailis} \runningtitle{Roughness of random paths by increment ratios}
\maketitle
\begin{abstract}
A general moment bound for sums of products of Gaussian
vector's functions extending the moment bound in Taqqu (1977, Lemma
4.5) is established. A general central limit theorem for triangular arrays of nonlinear functionals of multidimensional non-stationary
Gaussian sequences is proved. This theorem extends the previous results of Breuer and Major (1981), Arcones (1994) and others.
A Berry-Esseen-type bound in the above-mentioned central limit theorem is derived following Nourdin, Peccati and Podolskij (2011).
Two  applications of the above results are discussed. The first one refers to the asymptotic behavior of a roughness
statistic for  continuous-time Gaussian processes and the second one is
a central limit theorem satisfied by long memory locally stationary process.
\end{abstract}
\begin{quote}
\end{quote}

\vskip1cm

\section{Introduction}

This paper is devoted to the proof of two new results concerning functions of Gaussian vectors.
The first one (Lemma \ref{lemgauss} of Section \ref{moment}) is a moment bound for ``off-diagonal'' sums of
products of functions of Gaussian vectors in a general frame. It is an extension of an important lemma by Taqqu (1977, Lemma 4.5).
This result is useful for obtaining almost sure convergence and tightness of Gaussian subordinated functionals
and statistics, see Remark \ref{tight}  below.
The proof of Lemma \ref{lemgauss} uses the Hermitian decomposition of $\L^2$ function and the diagram formula.
A related but different moment bound is proved in
Soulier (2001, Corollary 2.1). % with the help of an elegant technique which does not use the diagram formula.
~\\

The second result is
a central limit
theorem (CLT) for arrays of random variables that are functions of Gaussian vectors, see Theorem \ref{tlcgauss} for a precise statement.
Theorem \ref{tlcgauss} generalizes and extends earlier results due to
Breuer and Major (1983), Giraitis and Surgailis (1985) and Arcones~(1994, Theorem 2)
to the case of non-stationary triangular arrays of Gaussian vectors. Extensions of the  Breuer-Major theorem were also obtained by Chambers and Slud (1989), Sanchez de Naranjo (1993) and
Nourdin {\it et al.} (2011).
%together with convergence rates.
Most of the above cited papers treat the case of a single stationary
Gaussian sequence and a function independent of $n$. Generalization to stationary or non-stationary triangular arrays
%Nourdin {\it et al.} (2011) improves the seminal Breuer and Major paper because this article also provides several upper bounds for Berry-Esseen type %relations (but only in the case of a stationary sequence of  Gaussian vectors).
%Such a result
is motivated by numerous statistical applications. Some examples of these applications, with a particular emphasis on strongly dependent
Gaussian processes,
are: statistics of time series (see for instance Bardet {\it et al.}, 2008,   Roueff and von Sachs, 2010), %nonparametric
kernel-type estimation of regression function (Guo and Koul, 2008),
nonparametric estimation of the local Hurst function of a continuous-time process %$(X_t)_{t\in [0,1]}$
from a discrete grid
$i/n, 1\le i \le n$
(Guyon and Leon, 1989, Bardet and Surgailis 2011, 2012).
Two particular applications (limit theorems for the Increment Ratio statistic of a Gaussian process admitting a tangent process and a CLT for functions of locally stationary Gaussian process) are discussed in Section \ref{Applications}. ~\\

Starting with the famous Lindeberg Theorem for independent random variables,
numerous studies devoted to
CLT for triangular arrays under various dependence conditions had appeared.
% was extensively studied
%under various dependence conditions.   %Proofs of such results are not only  straightforward applications of CLTs for sequences of random variables.
%Following the
%different assumptions on the dependence,
%We mention
%may cite as reference the book of
The case of martingale dependence was extensively studied in
Jacod and Shiryaev (1987).
%for martingale differences and the paper of
Rio (1995) discussed the case of strongly mixing sequences.
Some of more recent papers devoted to this question are Coulon-Prieur and Doukhan (2000)
(with a new weak dependence condition) and Dedecker and Merlev\`ede (2002) (with a necessary and sufficient condition for stable convergence of  normalized partial sums). The CLT for linear triangular arrays %({\it i.e.} sequence $(a_{in}\xi_i)_{1\leq i \leq n,~n\in \N}$
was discussed in detail in Peligrad and Utev (1997) for several forms of dependence conditions.
%on $(\xi_i)_i$. Note that all these results are obtained under a stationarity condition.
~\\

The case of Gaussian subordinated variables (functions of Gaussian vectors)
is rather exceptional among other dependence structures  since it allows for
very sharp conditions for CLT in terms of the decay rate of the covariance of Gaussian process and the Hermite rank of non-linear function. 
These conditions 
are close to being necessary and result in CLTs  ``in the vicinity" of non-central limit theorems, see Breuer and Major (1981), Arcones (1994),  
Dobrushin and Major (1979), Taqqu (1979).
The proofs of the above-mentioned results rely on
specific Gaussian techniques such as the Hermite expansion and the diagram formula; however, the recent paper
Nourdin {\it et al.} (2011) uses a different approach based on Malliavin's calculus and Stein's method, yielding also 
convergence rates in the CLT. The main 
difference between our Theorem \ref{tlcgauss} and the corresponding results in Arcones~(1994) and Nourdin {\it et al.} (2011) is that, contrary to these papers,
we do not assume stationarity of the underlying
Gaussian sequence  $\left({\mbf
Y}_{\! n}(k) \right)$ and discuss the case of subordinated sums $\sum_{k=1}^n f_{k,n}({\mbf Y}_{\! n}(k))$ where
$f_{k,n}$ may depend on $k$ and $n$. The last fact is important for statistical applications (see above). In the particular
case when $f_{k,n} = f$ do not depend on $k, n$ and $({\mbf Y}_{\! n}(k))$ is a stationary process independent of $n$,
Theorem \ref{tlcgauss} (iii) agrees with Arcones~(1994) and Nourdin {\it et al.} (2011, Theorem 1.1). The proof of Theorem \ref{tlcgauss}
%, see Section \ref{Section_CLT},
uses the diagram method and cumulants as in Giraitis and Surgailis (1985).  Section  \ref{breuer} obtains
a Berry-Esseen bound in this CLT using the approach and results in Nourdin {\it et al.} (2011). Let us note that a CLT 
for Gaussian subordinated arrays is also proved in Soulier (2001, Theorem 3.1); however, it requires
that Gaussian vectors are asymptotically independent and therefore his result is different from Theorem \ref{tlcgauss}. ~\\

\smallskip

\noindent {\bf Notation.} Everywhere below, ${\mbf X} = (X^{(1)}, \dots, X^{(\nu)})$ designates a
standardized Gaussian vector in $\R^\nu,  \, \nu \ge 1$, with zero mean $\E X^{(u)}  =0 $
and covariances $\E X^{(u)} X^{(v)} = \delta_{uv}, \,  u,v = 1,
\dots, \nu$. Letter $C$ stands for a constant whose precise value is unimportant and which may change from line
to line. The weak convergence of distributions is denoted by $\limiteloin $.

\section{A moment bound} \label{moment}

%Let ${\mbf X} = (X^{(1)}, \dots, X^{(\nu)})\in {\R}^\nu$ be a
%standardized Gaussian vector, with zero mean $\E X^{(u)}  =0 $
%and covariances $\E X^{(u)} X^{(v)} = \delta_{uv} \  (u,v = 1,
%\dots, \nu)$.
Let $\L^2 ({\mbf X}) $ denote the class of all
measurable functions $f = f({\mbf x}), {\mbf x} = (x^{(1)}, \dots,
x^{(\nu)}) \in {\R}^\nu$ such that $\|f\|^2 := \E f^2({\mbf X}) <
\infty $. For any multiindex ${\mbf k} = (k^{(1)}, \dots,
k^{(\nu)}) \in {\Z}^\nu_+ := \{(j^{(1)}, \dots, j^{(\nu)}) \in
{\Z}^\nu, j^{(u)} \ge 0 \ (1\le u \le \nu)\}$, let $H_{\mbf
k}({\mbf x}) = H_{k^{(1)}}(x^{(1)}) \cdots H_{k^{(\nu)}}(x^{(\nu)})$
be the (product) Hermite polynomial; $H_k(x) := (-1)^k
{\e}^{x^2/2}({\e}^{-x^2/2})^{(k)},$ $k=0,1, \ldots $ are standard
Hermite polynomials (with $({\e}^{-x^2/2})^{(k)}$ the $k$th derivative of the
function $x\mapsto {\e}^{-x^2/2}$).  Write  $|{\mbf k}| := k^{(1)}+ \dots +
k^{(\nu)},  \  {\mbf k}! := k^{(1)}! \cdots k^{(\nu)}!, \ {\mbf k} =
(k^{(1)}, \dots, k^{(\nu)}) \in {\Z}^\nu_+$. A function $f\in
\L^2({\mbf X})$ is said to have a Hermite rank $m\ge 0$ if $J_f
({\mbf k}) := \E f({\mbf X}) H_{\mbf k}({\mbf X}) =0 $  for any
${\mbf k}\in {\Z}_+^\nu, |{\mbf k}| < m$, and $J_f ({\mbf k}) \ne
0$ for some  ${\mbf k}, |{\mbf k}| = m$. It is well-known that any
$f \in \L^2({\mbf X}) $ having a Hermite rank $m\ge 0$ admits the
Hermite expansion
\begin{equation}\label{hermite0}
f({\mbf x}) = \sum_{|{\mbf k}|\ge m} \frac{J_f({\mbf k})}{{\mbf
k}!} H_{\mbf k}({\mbf x}),
\end{equation}
which converges in $\L^2 ({\mbf X}) $.\\
~\\
Let $({\mbf X}_1, \dots, {\mbf X}_n)$ be a collection of
standardized Gaussian vectors ${\mbf X}_t = (X^{(1)}_{t}, \dots,
X^{(\nu)}_{t}) \in {\R}^\nu $ having a joint Gaussian distribution
in ${\R}^{\nu n}$. Let $\varepsilon \in [0,1]$ be a fixed number.
Following Taqqu~(1977), we call $({\mbf X}_1, \dots, {\mbf X}_n)$
{\it $\varepsilon-$standard} if $|\E X^{(u)}_{t} X^{(v)}_{s}| \le
\varepsilon $ for any $t \neq s, 1\le t,s\le n$ and any $1\le u, v
\le \nu$.

\bigskip

As mentioned in the Introduction, Lemma \ref{lemgauss} generalizes 
Taqqu (1977, Lemma 4.5) to the case of a vector-valued
Gaussian family $({\mbf X}_1, \dots, {\mbf X}_{\! n}) $, taking
values in $ {\R}^\nu (\nu \ge 1)$. The lemma concerns the bound (\ref{momgauss}), below,  
%for $\sum\nolimits^\prime \big|\E \big[f_{1,t_1,n}({\mbf X}_{t_1}) \cdots
%f_{p,t_p,n}({\mbf X}_{t_p})\big] \big|$, 
where $f_{1,t,n}, \dots, f_{p,t,n}$
are square integrable functions among which the first $0\le \alpha
\le p$ functions $f_{1,t,n},\dots, f_{\alpha,t,n} $ for any $1\le
t \le n$ have a Hermite rank at least equal to $m\ge 1$ and where
$\sum\nolimits^\prime $ is the sum over all different indices
$1\le t_i \le n \ (1\le i \le p),  t_i \neq t_j (i\neq j)$. In the
case when $f_{j,t,n} = f_j $ does not depend on $t, n$, the bound (\ref{momgauss}) 
%Lemma \ref{lemgauss} 
coincides with that of Taqqu (1977, Lemma
4.5) provided $m \alpha $ is even, but  is worse than Taqqu's
bound in the more delicate case when $m \alpha  $ is odd. An advantage of
our proof is its relative simplicity (we do not use the
graph-theoretical argument as in Taqqu, 1977, but rather a simple
H\"older inequality). A different approach towards moment
inequalities for functions in vector-valued Gaussian variables is
discussed in Soulier (2001), leading to a different type of moment inequalities. 
% but is especially interesting when the Gaussian vectors are asymptotically independent and the components of the vectors are not independent.

\smallskip

\begin{lem}\label{lemgauss}
Let $({\mbf X}_1, \dots, {\mbf X}_n)$ be a $\varepsilon-$standard
Gaussian vector, ${\mbf X}_t = (X^{(1)}_{t}, \dots,X^{(\nu)}_{t})
\in {\R}^\nu, \, \nu \ge 1$, and let $f_{j,t,n} \in \L^2({\mbf X}),
1\le j \le p, \,  p\ge 2, \, 1\le t \le n $ be some functions. For
given integers $m\ge 1, 0\le \alpha \le p, n \ge 1$, define
\begin{eqnarray}
Q_n &:=&\max_{1\le t \le n} \sum_{1\le s\le n, s\neq t} \max_{1\le
u,v \le \nu} |\E X^{(u)}_{t} X^{(v)}_{s}|^m. \label{Q_N}
\end{eqnarray}
Assume that the functions $f_{1,t,n},\dots, f_{\alpha,t,n} $ have
a Hermite rank at least equal to $m$ for any $n\ge 1, 1\le t \le
n$, and that
\begin{equation}
\varepsilon < \frac{1}{\nu p-1}. \label{eps}
\end{equation}
Then
\begin{eqnarray}\label{momgauss}
\sum\nolimits^\prime \big|\E \big[f_{1,t_1,n}({\mbf X}_{t_1}) \cdots
f_{p,t_p,n}({\mbf X}_{t_p})\big] \big|
%\E|f_{1,t_1,n}({\mbf X}_{t_1}) \cdots
%f_{p,t_p,n}({\mbf X}_{t_p})| 
&\le& C( \varepsilon, p, m, \alpha,
\nu ) K n^{p- \frac{\alpha}{2}} Q_n^{\frac{\alpha}{2} },
\end{eqnarray}
where the constant $C( \varepsilon, p, m, \alpha, \nu ) $ depends
on $\varepsilon, p, m, \alpha, \nu $ only, and
\begin{equation}
K = \prod_{j=1}^p \max_{1\le t \le n}  \|f_{j,t,n}\|\quad \mbox{with}\quad \|f_{j,t,n}\|^2 = \E \big[f_{j,t,n}^2({\mbf X})\big].
\label{K_prod}
\end{equation}
%and the norms $ \|G_{j,t,N}\| = \E^{1/2} G^2_{j,t,N}({\mbf X}) \ (1\le j \le p, 1\le t \le N)$ only.
\end{lem}

\noindent {\it Proof.} Fix a collection $(t_1,\dots, t_p)$ of disjoint indices $t_i \neq t_j (i\neq j)$, and write $f_j = f_{j,
t_j, n},~ 1 \le j \le p$ for brevity.  Let $J_j({\mbf k}) :=
J_{f_j}({\mbf k}) =  \E \big[f_j({\mbf X}) H_{\mbf k}({\mbf X})\big] $ be
the coefficients of the Hermite expansion of $f_j$. Then,
\begin{eqnarray*}
|J_j({\mbf k})| &\le &\|f_j\| \prod_{i=1}^\nu \E^{1/2}
H^2_{k^{(i)}}(X) \\
&\le & \|f_j\| \prod_{i=1}^\nu (k^{(i)}!)^{1/2} = \|f_j\| ({\mbf
k}!)^{1/2}.
\end{eqnarray*}
Following Taqqu (1977, p. 213, bottom, p. 214, top), %for disjoint indices $t_i \neq t_j (i\neq j)$
we obtain
\begin{eqnarray*}
|\E f_1({\mbf X}_{t_1}) \cdots  f_p({\mbf X}_{t_p})|
&=&\left|\sum_{q=0}^\infty \sum_{|{\mbf k}_1|+ \dots +|{\mbf k}_p|
= 2q} \left\{\prod_{j=1}^p \frac{J_j({\mbf k}_j)}{{\mbf k}_j!}
\right\}
\E \big[H_{{\mbf k}_1}({\mbf X}_{t_1}) \cdots H_{{\mbf k}_p}({\mbf X}_{t_p})\big]\right| \\
&\le&K_1 \sum_{q=0}^\infty \sum_{|{\mbf k}_1|+ \dots +|{\mbf k}_p|
= 2q}
\frac{|\E H_{{\mbf k}_1}({\mbf X}_{t_1}) \cdots H_{{\mbf k}_p}({\mbf X}_{t_p})|}{ ({\mbf k}_1!\cdots {\mbf k}_p!)^{1/2}} \\
&\le&K_1 \sum_{q=0}^\infty \sum_{|{\mbf k}_1|+ \dots +|{\mbf k}_p|
= 2q}
\frac{\varepsilon^{(|{\mbf k}_1|+\dots + |{\mbf k}_p|)/2}
\E \prod_{1\le u \le \nu} \prod_{1\le j \le p} H_{k^{(u)}_{j}}(X)}{ ({\mbf k}_1!\cdots {\mbf k}_p!)^{1/2}}   \\
&\le&K_1 \sum_{q=0}^\infty \sum_{|{\mbf k}_1|+ \dots +|{\mbf k}_p|
= 2q} (\varepsilon (\nu p-1))^{(|{\mbf k}_1|+ \dots +|{\mbf
k}_p|)/2} < \infty,
\end{eqnarray*}
where $X \sim {\mathcal N}(0,1)$ and
$$
K_1 :=   \|f_{1,t_1,n}\| \cdots \|f_{p,t_p,n}\| \ \le \  K,
$$
where $K$ is defined in (\ref{K_prod}) and $K$ is independent of
$t_1, \dots, t_p$,  and where we used the assumption  (\ref{eps})
to get the convergence of the last series. Therefore,
\begin{eqnarray*}
\sum\nolimits^\prime \big| \E \big[f_{1,t_1,n}({\mbf X}_{t_1}) \cdots
f_{p,t_p,n}({\mbf X}_{t_p})\big] \big| &\le& K \sum_{q=0}^\infty
\sum_{\tiny
\begin{array}{c}
|{\mbf k}_1|+ \dots +|{\mbf k}_p| = 2q\\
|{\mbf k}_1| \ge m, \dots,  |{\mbf k}_\alpha| \ge m
\end{array}}
\sum\nolimits^\prime \frac{|\E H_{{\mbf k}_1}({\mbf X}_{t_1}) \cdots
H_{{\mbf k}_p}({\mbf X}_{t_p})| } {({\mbf k}_1!\cdots {\mbf
k}_p!)^{1/2}}.
\end{eqnarray*}
Now, the following bound remains to be proved: for any integers
$m\ge 1, 0\le \alpha \le p, n \ge 1$ and any multiindices ${\mbf
k}_1, \dots {\mbf k}_p \in {\Z}^\nu_+ $ satisfying $ |{\mbf k}_1|+
\dots +|{\mbf k}_p| = 2q, \ |{\mbf k}_1| \ge m, \dots, |{\mbf
k}_\alpha| \ge m $,
\begin{equation}
\sum\nolimits^\prime |\E H_{{\mbf k}_1}({\mbf X}_{t_1}) \cdots
H_{{\mbf k}_p}({\mbf X}_{t_p})| \ \le \  C_1 (\varepsilon (\nu
p-1))^{(|{\mbf k}_1|+ \dots +|{\mbf k}_p|)/2} ({\mbf k}_1!\cdots {\mbf
k}_p!)^{1/2} n^{p-\frac{\alpha}{2}}Q_n^{\frac{\alpha}{2} },
\label{hermite_bdd}
\end{equation}
where $C_1 $ is some constant depending only on $p, \nu, \alpha,
\varepsilon$, and independent of ${\mbf k}_1, \dots, {\mbf k}_p, n$.

\smallskip

First, we write the expectation on the left hand side of
(\ref{hermite_bdd}) as a sum of contributions of diagrams. Let
\begin{equation}
T:= \pmatrix{ (1,1) &(1,2) &\dots &(1,k_1)\cr (2,1) &(2,2) &\dots
&(1,k_2)\cr \dots \cr (p,1) &(p,2) &\dots &(p,k_p)\cr}
\label{T_table}
\end{equation}
be a table having $p$ rows $\tau_1, \dots, \tau_p$ of respective
lengths $|\tau_u| = k_u = |{\mbf k}_u| = k^{(1)}_{u}+ \dots +
k^{(\nu)}_{u}$ (we write $T = \bigcup_{u=1}^p \tau_u)$. A {\it
sub-table} of $T$ is a table $T' = \bigcup_{u \in U} \tau_u, \  U
\subset \{1, \dots, p\} $ consisting of some rows of $T$ written
from top to bottom in the same order as  rows in $T$; clearly any
sub-table $T'$ of $T$ can be identified with a (nonempty) subset $U
\subset \{1, \dots, p\}$. A {\it diagram} is a partition $\gamma $
of the table $T$ by pairs (called {\it edges} of the diagram) such
that no pair belongs to the same row. A diagram $\gamma $ is
called {\it connected} if the table $T$  cannot be written as a
union $T =T' \cup T''$ of two disjoint sub-tables $T', T''$ so that
$T' $ and $T''$ are partitioned by $\gamma $ separately. Write
$\Gamma (T), \Gamma_c(T)$ for the class of all diagrams and the
class of all connected diagrams over the table $T$, respectively.
Let
$$
\rho (t,s) := \max_{1\le u,v \le \nu} |\E X^{(u)}_{t} X^{(v)}_{s}|
\qquad (t\ne s).
$$
Note $ 0\le \rho (t,s) \le \varepsilon $ and $Q_n =\max_{1\le t \le
n} \sum_{1\le s\le n, s\neq t} \rho^m(t,s). $ By the diagram formula
for moments of Hermite (Wick) polynomials (see e.g. Surgailis,
2000),
\begin{eqnarray}
|\E H_{{\mbf k}_1}({\mbf X}_{t_1}) \cdots H_{{\mbf k}_p}({\mbf
X}_{t_p})|
&\le&\sum_{\gamma \in \Gamma (T)} \prod_{1 \le u <v \le p} (\rho (t_u,t_v))^{\ell_{uv}} \label{T} \\
&=&\sum_{(U_1, \dots, U_h)} \prod_{r=1}^h  \sum_{\gamma \in
\Gamma_c(U_r)} \prod_{u,v \in U_r, u<v}  (\rho
(t_u,t_v))^{\ell_{uv}}, \label{U}
\end{eqnarray}
where $\ell_{uv} $ is the number of edges between rows $\tau_u$
and $\tau_v $ in the diagram $\gamma $ over  table $T$, and the
sum $\sum_{(U_1, \dots, U_h)}$ is taken over all partitions $(U_1,
\dots, U_h), h = 1, 2, \dots, [p/2] $ of $\{1, \dots, p\} $ by
nonempty subsets $U_r$ of cardinality $|U_r| \ge 2$. (Thus,
(\ref{U}) follows from (\ref{T}) by decomposing $\gamma \in \Gamma
(T)$ into connected components $\gamma_r \in \Gamma_c(U_r), r=1,
\dots, h; h =1,\dots, [p/2]$; the restriction $|U_r|\ge 2 $ stems
from the fact that any edge must necessarily connect different
rows.)   From (\ref{U}) we obtain
\begin{equation}
\sum\nolimits^\prime \big|\E \big[f_1({\mbf X}_{t_1}) \cdots  f_p({\mbf
X}_{t_p})\big] \big|  \ \le \  \sum_{(U_1, \dots, U_h)} \prod_{r=1}^h
\sum_{\gamma \in \Gamma_c(U_r)} I_{n,U_r}(\gamma), \label{G_I}
\end{equation}
where, for any sub-table $U \subset T $ having at least two rows
and for any connected diagram $\gamma \in \Gamma_c(U)$,  the
quantity $I_{n,U}(\gamma)$ is defined by
$$
I_{n,U}(\gamma) :=    \sum\nolimits^\prime \prod_{u,v \in U, u<v}
(\rho (t_u,t_v))^{\ell_{uv}}
$$
where (recall) the product is taken over all ordered pairs of rows
$(\tau_u, \tau_v), u<v $ of the table $U$, and $\ell_{uv}$ is the
number of edges in $\gamma $ between the $u$th and the $v$th rows.
Below we  prove the bound
\begin{equation}
I_{n,U}(\gamma) \le K_3 \epsilon^{|{\mbf
k}_U|/2}n^{|U|-\frac{\alpha(U)}{2}} (n Q_n)^{\frac{\alpha(U)}{2}},
\label{I_U}
\end{equation}
where $|{\mbf k}_U| := \sum_{u\in U} k_u $ is the number of points
of table $U$ and  $\alpha(U) := |\{1, \dots, \alpha\} \bigcap U| =
\# \{u\in U: |{\mbf k}_u| \ge m \} $ is the number of rows in $U$
having at least $m$ points. Clearly, it suffices to show
(\ref{I_U}) for $U=T$.

Next, let for $1\le u,v \le p, u \ne v $, denote
\begin{eqnarray}
R_{uv} &:=&\left(\sum\nolimits_{1\le t\le n}
\left(\sum\nolimits_{1\le s \le n, s\neq t} \rho^{k_u}(s,t)
\right)^{k_v/k_u} \right)^{\ell_{uv}/k_v}.  \label{R_uv}
\end{eqnarray}
Let $A := \{1, \dots, \alpha\}, A' := \{1, \dots,p\}\backslash A =
\{\alpha+1, \dots, p\} $. It follows immediately from the definition of
$R_{uv}$ and $\rho(s,t) $ that
\begin{eqnarray}
R_{uv}&\le&\cases{n^{\frac{\ell_{uv}}{k_v}}
Q_n^{\frac{\ell_{uv}}{k_u}}
\varepsilon^{(1-\frac{m}{k_u})\ell_{uv}}, &if $u\in A$, \cr
n^{\frac{\ell_{uv}}{k_u} + \frac{\ell_{uv}}{k_v}}
\varepsilon^{\ell_{uv}}, &if $u \in A^c$. \cr} \label{R_bdd}
\end{eqnarray}
By the H\"older inequality (see Giraitis and Surgailis, 1985, p.202,
for details),
\begin{eqnarray}
I_{n,T}(\gamma) &\le& \min\left(\prod_{1\le u<v\le p} R_{uv},
\prod_{1\le u<v\le p} R_{vu} \right). \label{H_bdd}
\end{eqnarray}
For any subset $U\subset \{1, \dots, p\}$, let
\begin{equation}
L(U) := \sum_{u\in U}\sum_{u<v\le p} \frac{\ell_{uv}}{k_u}, \qquad
L^*(U) := \sum_{u\in U} \sum_{1\le v <u} \frac{\ell_{uv}}{k_u},
\label{L(T)}
\end{equation}
$L := L(T), L^* := L^*(T)$. Clearly,
\begin{equation}
L(U) + L^*(U) \  = \  \sum_{u\in U} \frac{1}{k_u} \sum_{v=1,
\dots, p, v \ne u} \ell_{uv} \  = \  |U|  \label{U_iden}
\end{equation}
is the number of points in $U$. From (\ref{R_bdd}) -
(\ref{H_bdd}),
\begin{eqnarray*}
I_{n,T}(\gamma) \ \le \  \min\left(n^{L^* + L(A^c)} Q_n^{L(A)}
\varepsilon^{|T|/2 - mL(A)}, n^{L + L^*(A^c)} Q_n^{L^*(A)}
\varepsilon^{|T|/2 - mL^*(A)}\right),
\end{eqnarray*}
where $|T| = \sum_{u=1}^p k_i$. As $0\le L(A),L^*(A) \le p $, see
(\ref{U_iden}), we obtain
\begin{eqnarray*}
I_{n,T}(\gamma)
&\le&%\varepsilon^{|T|/2 - mp} \min \left(N^{L^*  + L(A^c)} Q_N^{L(A)}, N^{L + L^*(A^c)} Q_N^{L^*(A)} \right) \\
\varepsilon^{|T|/2 - mp}\min \left(n^{L^*(A) + L^*(A^c) +
L(A^c)} Q_n^{L(A)}, n^{L(A) + L(A^c) + L^*(A^c)} Q_n^{L^*(A)} \right) \\
&=&\varepsilon^{|T|/2 - mp}n^{p-\alpha} \min \left(n^{L^*(A)} Q_n^{L(A)}, n^{L(A)} Q_n^{L^*(A)} \right) \\
&=&\varepsilon^{|T|/2 - mp} n^{p-\alpha} (n
Q_n)^{\frac{\alpha}{2}} \min
\left(\left(n/Q_n\right)^{\frac{\alpha}{2}- L(A)},
\left(n/Q_n\right)^{L(A)-\frac{\alpha}{2}}
\right) \\
&\le &\varepsilon^{|T|/2 - mp}n^{p-\alpha} (n
Q_n)^{\frac{\alpha}{2}},
\end{eqnarray*}
proving (\ref{I_U}).

With (\ref{I_U})-(\ref{G_I}) in mind,
\begin{eqnarray*}
\sum\nolimits^\prime |\E H_{{\mbf k}_1}({\mbf X}_{t_1}) \cdots
H_{{\mbf k}_p}({\mbf X}_{t_p})| &\le&
C_3 \varepsilon^{|T|/2} \sum_{(U_1, \dots, U_h)} \prod_{r=1}^h
\sum_{\gamma \in \Gamma_c(U_r)} n^{|U_r| - \frac{\alpha (U_r)}{2}} Q_n^{\frac{\alpha (U_r)}{2}}\\
&=& C_3 \varepsilon^{|T|/2}n^{p- \frac{\alpha}{2}}
Q_n^{\frac{\alpha}{2}}
\sum_{(U_1, \dots, U_h)} \prod_{r=1}^h \sum_{\gamma \in \Gamma_c(U_r)}1 \\
&=& C_3 \varepsilon^{|T|/2}n^{p- \frac{\alpha}{2}}
Q_n^{\frac{\alpha}{2}}  \sum_{\gamma \in \Gamma (T)}1,
\end{eqnarray*}
where the last sum (= the number of all diagrams over the table
$T$) does not exceed
$$
|\E H_{k^{(1)}_{1}}(X) \cdots H_{k^{(\nu)}_{1}}(X) \cdots
H_{k^{(1)}_{p}}(X) \cdots H_{k^{(\nu)}_{p}}(X)| \le (p\nu -1)^{(|{\mbf
k}_1|+\dots + |{\mbf k}_p|)/2} ({\mbf k}_1! \cdots {\mbf k}_p!)^{1/2},
$$
see  Taqqu~(1977, Lemma 3.1). This proves the bound
(\ref{hermite_bdd}) and the lemma, too.  \hfill $\Box$ \\

Lemma \ref{lemgauss} can be extended to non-standardized Gaussians as follows.  To this end, we introduce some
definitions. Let ${\mbf Y} = (Y^{(1)}, \dots, Y^{(\nu)}) \in \R^\nu $ be a Gaussian vector with zero mean and
non-degenerate covariance matrix $\Sigma = (\E Y^{(u)} Y^{(v)})_{1\le u,v \le \nu}$.
Let $\L^2 ({\mbf Y}) $ denote the class of all
measurable functions $f: \R^\nu \to \R$
%   = G({\mbf x}), {\mbf x} = (x^{(1)}, \dots,
%x^{(\nu)}) \in {\R}^\nu$
with $%\|G\|^2_\Sigma :=
\E f^2({\mbf Y}) <
\infty $. Let $m \ge 0$ be an integer.
We say that {\it  $f \in \L^2 ({\mbf Y})$ has a generalized Hermite rank not less than $m$} if either
$m = 0,$  or $m \ge 1 $ and
\begin{equation}\label{PG0}
\E [ P({\mbf Y}) f({\mbf Y}) ]  = 0 \qquad \text{for all} \  P \in {\cal P}_{m-1}({\R}^\nu)
\end{equation}
hold, where ${\cal P}_m({\R}^\nu)$ stands for  the class of all polynomials $P$ in variables $y^{(1)}, \dots, y^{(\nu)}$ of degree $m$,
that is,
$P({\mbf y}) = \sum_{0\le |{\mbf j}|\le m} c({\mbf j}) {\mbf y}^{\mbf j}
= \sum_{j^{(1)} \ge 0, \dots, j^{(\nu)} \ge 0: j^{(1)} + \dots + j^{(\nu)} \le m} $  $ c(j^{(1)}, \dots, j^{(\nu)}) $
$(y^{(1)})^{j^{(1)}} \dots (y^{(\nu)})^{j^{(\nu)}} $.

Let ${\mbf
X} :=  \Sigma^{-1/2} {\mbf
Y}, \,   \tilde f({\mbf x}) := f( \Sigma^{1/2} {\mbf x}).    $ Then   ${\mbf
X}$ has a standard Gaussian distribution in $\R^\nu$ and $\tilde f \in \L^2 ({\mbf X})$
with
\begin{equation}\label{PG2}
\|\tilde f \|^2 =  \E|\tilde f({\mbf X})|^2 =  \E|f ({\mbf Y})|^2.
\end{equation}

The following proposition is known, see Nourdin {\it et al.} (2011, Proposition 2.1),
 Soulier (2001, p.195), but we include a proof of it for completeness.

\begin{prop} \label{GHrank} Let ${\mbf Y},  {\mbf X}, \, f \in \L^2 ({\mbf Y}), \, \tilde f \in \L^2 ({\mbf X})$ be defined  as above and
$m \ge 0$ be a given integer. $f $ has a generalized Hermite rank not less than $m$ if and only if
$\tilde f$ has a Hermite rank not less than $m$.

\end{prop}

\noindent {\it Proof.} The above proposition is true if ${\mbf Y} =  {\mbf X}$ has a standard Gaussian distribution;
see Soulier (2001, p.194). By definition
 %(\ref{PG0}) can be rewritten as
\begin{equation}\label{PG3}
\E [ P({\mbf Y}) f({\mbf Y}) ]  =  \E [ \tilde P({\mbf X}) \tilde f({\mbf X}) ],
\end{equation}
where
$\tilde P({\mbf x}) :=  P(\Sigma^{1/2} {\mbf x}) $.  Clearly,  $P  \in {\cal P}_{m-1}(\R^\nu)$ implies that
$\tilde  P   \in {\cal P}_{m-1}(\R^\nu)$ is a polynomial of degree $m-1$. Therefore
$\tilde f$ having a Hermite rank not less than $m$ implies by (\ref{PG3}) that  $f$ has a generalized Hermite rank not less than $m$.
The converse statement again follows from (\ref{PG3}), by taking $P({\mbf y}) = \hat P(\Sigma^{-1/2} {\mbf y})$, where
$\hat P \in {\cal P}_{m-1}(\R^\nu)$ is an arbitrary polynomial of degree $m-1$. \hfill $\Box$ \\
%Proposition \ref{GHrank} is proved.

%\medskip

Let $({\mbf Y}_{\!1}, \dots, {\mbf Y}_{\!n})$ be a collection of
Gaussian vectors ${\mbf Y}_{\!t} = (Y^{(1)}_{t}, \dots,
Y^{(\nu)}_{t}) \in {\R}^{\nu} $ with zero mean $\E {\mbf Y}_{\!t} = 0 $ and non-degenerated covariance matrices
$\Sigma_t = \big({\rm Cov}\big(Y^{(u)}_{t}, Y^{(v)}_{t}\big)\big)_{1\le u,v\le \nu}, $
having a joint Gaussian distribution
in ${\R}^{\nu n}$. Let $\varepsilon \in [0,1]$ be a fixed number.
Call $({\mbf Y}_{\!1}, \dots, {\mbf Y}_{\!n})$
{\it $\varepsilon-$correlated} if $\big|{\rm Cor}\big(Y^{(u)}_{t}, Y^{(v)}_{s}\big)\big| \le
\varepsilon $ for any $t \neq s, 1\le t,s\le n$ and any $1\le u, v
\le \nu$. Clearly, if the ${\mbf Y}_{\!t}$'s are standard, this is equivalent to
 $({\mbf Y}_{\!1}, \dots, {\mbf Y}_{\!n})$ being
$\varepsilon-$standard.

We also use some elementary facts about matrix norms. Let $|{\mbf x}| = \big(\sum_{i=1}^\nu (x^{(i)})^2 \big)^{1/2} $ denote the Euclidean norm
in $\R^\nu$,
$A = (a_{ij}) $  a real $\nu\times \nu-$matrix, $A^\intercal $ the transposed matrix,
$I$ the unit matrix,
and
$\|A \| := \sup_{|{\mbf x}| =1} |A {\mbf x}|$  the matrix spectral norm, respectively.
 Then $\|A\|_{\infty} :=   \max_{1\le i,j\le \nu} |a_{ij}| \le \|A\| \le \nu \|A\|_{\infty}$ and
  $\|A B \| \le \|A \| \,\|B \|$ for any such matrices
$A, B$. An orthogonal matrix  $O = (o_{ij}) $ satisfies $O O^\intercal =  O^\intercal O = I $ and $\| O\| = \| O^\intercal\| =  1.$
Any symmetric matrix $A$ can be written as $A = O^\intercal \Lambda O,$ where $O $ is an orthogonal matrix and
$\Lambda $ is a diagonal matrix. In addition, if $A$ is positive definite, then $\| A\| = \| \Lambda \| =
\lambda_{\max}, \,
\| A^{-1} \| = \| \Lambda^{-1} \| = \lambda^{-1}_{\min}, $ where $\lambda_{\max} \ge \lambda_{\min} >0$ are the largest
and the smallest eigenvalues of $A$. We shall also  use the facts that for
any symmetric positive definite matrix $A$,
\begin{equation}\label{matrixroot}
\|A^{1/2}\|  = \| A \|^{1/2}, \qquad   \|A^{-1/2}\|  = \| A^{-1} \|^{1/2},
\end{equation}
since $\|A^{1/2}\| = \|O^\intercal \Lambda^{1/2} O\|
= \| \Lambda^{1/2} \| =
\lambda^{1/2}_{\max}, \,  \|A^{-1/2}\| = \|O^\intercal \Lambda^{-1/2} O\|
= \| \Lambda^{-1/2} \| =
\lambda^{-1/2}_{\min}$.
%{\bf [Are these facts too well known ? - D.]}

\smallskip

\begin{cor}\label{corgauss}
Let $({\mbf Y}_1, \dots, {\mbf Y}_n)$ be an $\varepsilon-$correlated
Gaussian vector, ${\mbf Y}_t = (Y^{(1)}_{t}, \dots, Y^{(\nu)}_{t})
\in {\R}^\nu \ (\nu \ge 1)$, with zero mean $\E {\mbf Y}_{\!t} = 0 $ and non-degenerated covariance matrices
$\Sigma_t $ satisfying
\begin{equation}\label{cmax}
\max_{1 \le t \le n} %\big( \| \Sigma_t \| +
\| \Sigma^{-1}_t \|  \ \le \ c_{\max}
\end{equation}
for some constant $c_{\max} >0$.
Let $f_{j,t,n} \in \L^2({\mbf Y}_t),
1\le j \le p \  (p\ge 2), 1\le t \le n $ be some functions. For
given integers $m\ge 1,
0\le \alpha \le p, n \ge 1$,
let $Q_n$ denote the sum in (\ref{Q_N}) where $ X^{(u)}_{t}, \, X^{(v)}_{s} $ are replaced
by $Y^{(u)}_{t}, Y^{(v)}_{s}$, respectively.
%\begin{eqnarray}
%Q_N &:=&\max_{1\le t \le N} \sum_{1\le s\le N, s\neq t} \max_{1\le
%u,v \le \nu} |\E Y^{(u)}_{t} Y^{(v)}_{s}|^m. \label{Q_NY}
%\end{eqnarray}
Assume that the functions $f_{1,t,n},\dots, f_{\alpha,t,n} $ have
a generalized Hermite rank at least equal to $m$ for any $n\ge 1, 1\le t \le
n$, and that
\begin{equation}
\varepsilon   < \frac{1}{(\nu p-1)\nu^2 c_{\max} }. \label{eps1}
\end{equation}
Then
\begin{eqnarray*}
\sum\nolimits^\prime \big| \E \big[ f_{1,t_1,n}({\mbf Y}_{t_1}) \cdots
f_{p,t_p,n}({\mbf Y}_{t_p})\big] \big| &\le& C K n^{p- \frac{\alpha}{2}} Q_n^{\frac{\alpha}{2} },
\end{eqnarray*}
where $K := \prod_{j=1}^p \max_{1\le t \le n} \E^{1/2} \big[ f_{j,t,n}^2({\mbf Y}_t)\big] $ and the constant $ C = C( \varepsilon, p, m, \alpha, \nu, c_{\max} ) $ depends
on $\varepsilon, p, m, \alpha, $  $\nu, c_{\max} $ only.

%and
%\begin{equation1}
%$K = \prod_{j=1}^p \max_{1\le t \le N} \E^{1/2} G_{j,t,N}^2({\mbf Y}_t). $
%\label{K_prod}
%\end{equation}
%and the norms $ \|G_{j,t,N}\| = \E^{1/2} G^2_{j,t,N}({\mbf X}) \ (1\le j \le p, 1\le t \le N)$ only.
\end{cor}

\noindent {\it Proof.} We will reduce the above inequality to that of Lemma \ref{lemgauss}, as follows. Let
${\mbf X}_t := \Sigma^{-1/2}_t {\mbf Y}_t, \, \tilde f_{j,t,n}({\mbf x}) :=   f_{j,t,n}(\Sigma^{1/2}_t {\mbf x})$. The
${\mbf X}_t$'s have a standard Gaussian distribution in $\R^\nu$
and the $\tilde f_{j,t,n}$'s satisfy $\|\tilde  f_{j,t,n}\|^2
= \E \big[f_{j,t,n}^2({\mbf Y}_t)\big]$, see (\ref{PG2}). By Proposition \ref{GHrank},
$\tilde f_{j,t,n}, j=1, \dots, \alpha $ have a Hermite rank not less than $m$.
Next, using (\ref{matrixroot}), (\ref{cmax})
and the fact that the ${\mbf Y}_{\!t}$'s are $\varepsilon-$correlated, for any
$t\ne s, \, 1 \le t, s \le n, 1\le u, v \le \nu $
\begin{equation}\label{XYbdd}
|\E X^{(u)}_t X^{(v)}_s |\ \le \ \nu^2 \|\Sigma^{-1/2}_t \|_\infty \|\Sigma_s^{-1/2} \|_\infty \max_{1\le u, v \le \nu}
|\E Y^{(u)}_t Y^{(v)}_s | \ \le  \ \varepsilon \nu^2 \|\Sigma^{-1/2}_t \| \|\Sigma_s^{-1/2} \|
\ \le  \  \varepsilon \nu^2 c_{\max}
\end{equation}
This  implies that the Gaussian vector $({\mbf X}_{\!1}, \dots, {\mbf X}_{\!n})\in {\R}^{\nu n}$ is $\tilde \varepsilon-$standard,
where  $ \tilde \varepsilon := \varepsilon \nu^2 c_{\max}$. Then, in view of (\ref{eps1}),
(\ref{momgauss}) of Lemma \ref{lemgauss} applies, according to which
\begin{eqnarray*}
&&\sum\nolimits^\prime \big| \E \big[f_{1,t_1,n}({\mbf Y}_{t_1}) \cdots
f_{p,t_p,n}({\mbf Y}_{t_p})\big] \big|\ =\  \sum\nolimits^\prime \big| \E \big[\tilde f_{1,t_1,n}({\mbf X}_{t_1}) \cdots
\tilde f_{p,t_p,n}({\mbf X}_{t_p})\big] \big| \\
&& \le\  C(\tilde \varepsilon, p, m, \alpha,
\nu )\tilde K n^{p- \frac{\alpha}{2}}\tilde  Q_n^{\frac{\alpha}{2} } \ \le  \
C(\varepsilon, p, m, \alpha,
\nu, c_{\max} ) K n^{p- \frac{\alpha}{2}} Q_n^{\frac{\alpha}{2} },
\end{eqnarray*}
where $\tilde K, \tilde Q_n$ are the corresponding quantities in  Lemma \ref{lemgauss}  (\ref{momgauss}) satisfying
$\tilde K = K, \tilde Q_n \le  (\varepsilon \nu^2 c_{\max})^m Q_n$ by (\ref{PG2}), (\ref{XYbdd}), respectively.
\hfill $\Box$ \\
%Corollary \ref{corgauss} is proved.

%{\bf [It seems the above proof used only $\max_{1\le t \le N} \| \Sigma^{-1}_t \| < \ c_{\max} $ instead of
%(\ref{cmax}), which is a bit strange - the matrices $\Sigma_t$ need not be bounded? Note also that
%(\ref{eps1}) is not optimal since it does not reduce to (\ref{eps}) in the $\epsilon-$standard case. The last fact could be due
%to (\ref{XYbdd}) being not optimal
%- D.] }

We remark that condition (\ref{eps1}) is not optimal since it does not reduce to (\ref{eps}) in the $\varepsilon-$standard case.
This loss of optimality is due to the use of robust inequalities for matrix norms in (\ref{XYbdd}).

\begin{rem} \label{tight} {\rm As mentioned in the Introduction, Lemma \ref{lemgauss} and Corollary \ref{corgauss} can be used
for proving the tightness and the strong law of large numbers of various non-linear
statistics from Gaussian observations. See Bardet and Surgailis (2011, 2012) on application for
roughness estimation and Cs\"org{\H o} and Mielnichuk (1996), Koul and Surgailis (2002) for empirical process. The above-mentioned 
applications concern the 4th moment bound $M_{n} := \E \big(\sum_{t=1}^n f_{t,n}({\mbf Y}_{\! n}(t)) \big)^4 = O(n^{-\kappa})$
for a suitable $\kappa >0$, where $({\mbf Y}_{\! n}(t)), (f_{t,n})$ satisfy similar conditions as in
Corollary \ref{corgauss}. Clearly, $M_{n} = \sum_{t_1, \dots, t_4=1}^n   \E \big[\prod_{i=1}^4 f_{t_i,n}({\mbf Y}_{\! n}(t_i)) \big] $ can be decomposed into four terms according to the number of coinciding
``diagonals'' $t_i = t_j $ in the last sum, where each term can be estimated  with the help of Corollary \ref{corgauss}. 
Let us note that condition (\ref{eps1}) in the above applications is guaranteed by a preliminary ``decimation'' of the sum $\sum_{t=1}^n f_{t,n}({\mbf Y}_{\! n}(t),$ see 
(Cs\"org{\H o} and Mielnichuk, 1996) and (Bardet and Surgailis, 2012)  for details.

% \le C \big(\sum_{t_1,t_2, t_3, t_4}^{\prime} +  \sum_{t_1,t_2, t_3}^{\prime} +
%\sum_{t_1,t_2}^{\prime} + \sum_{t_1}^\prime\big)$

}
\end{rem}

\section{A CLT for triangular array of functions of Gaussian vectors} \label{Section_CLT}
Let $\left({\mbf
X}_{\! n}(k) \right)_{1\le k \le n, n \in {\N}} $ be a  triangular
array of standardized Gaussian vectors with values in ${\R}^\nu, \
{\mbf X}_{\! n}(k) = (X^{(1)}_n(k), \dots, X^{(\nu)}_n(k)), \  \E
X^{(p)}_n(k) = 0, \  \E X^{(p)}_n(k) X^{(q)}_n(k) = \delta_{pq}$.
Now define,
$$
r^{(p,q)}_n(j,k) :=  \E X^{(p)}_n(j) X^{(q)}_n(k)  \qquad (1 \le
j, k \le n).
$$
For a given integer $m\ge 1$, introduce the following assumptions:
for any $1\le p, q \le \nu$,
\begin{eqnarray}
\sup_{n \ge 1} ~ \max_{1\le k \le n} \sum_{1 \leq j \leq n} \big|
r^{(p,q)}_n(j,k)\big|^m \hspace{-5mm} &<& \hspace{-5mm} \infty,
\label{cov_max} \\
\sup_{n \ge 1} ~ \frac 1 n \hspace{-0.5cm} \sum_{\tiny
\begin{array}{c}
 1\le j,k \le n
\\
 |j-k| >K\end{array}} \left|r^{(p,q)}_n(j,k) \right|^m \hspace{-5mm} &\limiteK& \hspace{-5mm} 0,
 \label{cov_rank_K} \\
\forall (j, k)\in \{1,\ldots,n\}^2,~~~~\left|
r^{(p,q)}_n(j,k)\right| \hspace{-5mm} &\le&
\hspace{-5mm}|\rho(j-k)|~~~~\mbox{with} ~~\sum_{j\in {\Z}}
|\rho(j)|^m < \infty. \label{rho_dom}
\end{eqnarray}
Note (\ref{rho_dom}) $\Rightarrow$ (\ref{cov_max}) and (\ref{rho_dom})
$\Rightarrow$ (\ref{cov_rank_K}). Let $\L^2_0({\mbf X}) := \{f
\in \L^2({\mbf X}): \E f({\mbf X}) = 0\}, $ where ${\mbf X}\in \R^\nu$ denotes a standard Gaussian vector as above.

\begin{thm}\label{tlcgauss}
Let $\left({\mbf X}_{\! n}(k) \right)_{1\le k \le n, n \in {\N}} $
be a  triangular array of standardized Gaussian
vectors.
\begin{enumerate}
\item [(i)] Assume (\ref{cov_max}).
Let $f_k\in \L^2_0({\mbf X}) \ (1\le k \le n)$ have a Hermite rank at
least $m \in \N^*$.
%with  $\L^2_0({\mbf X}) = \{f
%\in \L^2({\mbf X}): \E f({\mbf X}) = 0\}$.
Then there exists a constant $C$ independent of
$n$ and $f_k, 1\le k \le n$ such that
\begin{equation}
\E \Big(n^{-1/2} \sum_{k=1}^n f_k\left({\mbf X}_{\! n}(k)\right)
\Big)^2 \  \leq \ C \max_{1\le k \le n} \|f_k\|^2.
\label{f_bdd}
\end{equation}
\item [(ii)] Assume (\ref{cov_max}) and (\ref{cov_rank_K}).
Let $f_{k,n} \in  \L^2_0({\mbf X}) \ (n \ge 1, 1 \le k \le n)$ be a
triangular array of functions all having Hermite rank at least $m
\in \N^*$. Assume that there exists a $ \L^2_0({\mbf X})-$valued
continuous function $\phi_\tau, \tau \in [0,1]$, such that
\begin{equation}
\sup_{\tau \in (0,1]} \|f_{[\tau n],n} - \phi_{\tau}\|^2 =
\sup_{\tau \in (0,1]} \E \big(f_{[\tau n],n}({\mbf X}) -
\phi_{\tau}({\mbf X})\big)^2 \limiten 0. \label{f_n}
\end{equation}
Moreover, let
\begin{equation}
\sigma^2_n := \E \Big(n^{-1/2} \sum_{k=1}^n f_{k,n}\left({\mbf
X}_{\! n}(k)\right) \Big)^2\limiten \sigma^2, \label{sigma_n}
\end{equation}
where $\sigma^2 >0$. Then
\begin{equation}
n^{-1/2} \sum_{k=1}^n f_{k,n}\left({\mbf X}_{\! n}(k)\right)
\limiteloin {\mathcal N}(0, \sigma^2). \label{CLT}
\end{equation}
\item [(iii)] Assume
(\ref{rho_dom}). Moreover, assume that for any $\tau \in [0,1]$ and
any $J \in \N^*$,
\begin{equation}
\big({\mbf X}_{\! n}([n\tau]+j)\big)_{-J \le j \le J} \
\limiteloin \  \big({\mbf W}_{\! \tau}(j)\big)_{-J \le j \le J}
,   \label{Y_W}
\end{equation}
where $ \left({\mbf W}_{\! \tau}(j)\right)_{j\in {\Z}}$ is a
stationary Gaussian process taking values in ${\R}^\nu$ and
depending on parameter $\tau \in (0,1)$. Let  $f_{k,n} \in
\L^2_0({\mbf X}) \ (n \ge 1, 1 \le k \le n)$ satisfy the same
conditions as in part (ii), with exception of (\ref{sigma_n}).
Then (\ref{sigma_n}) and (\ref{CLT}) hold, with
\begin{equation}
\sigma^2 \ = \  \int_0^1 \d \tau  \Big(\sum_{j\in {\Z}} \E \left[
\phi_\tau \left({\mbf W}_{\! \tau}(0)\right) \phi_\tau \left({\mbf
W}_{\! \tau}(j)\right)\right ]\Big).  \label{sigma}
\end{equation}
\end{enumerate}
\end{thm}
%In Bardet and Surgailis (2008 and 2009) applications of this theorem to

We remark that parts (i) and (ii) of Theorem \ref{tlcgauss} are natural extensions of Theorem 2 of Arcones (1994) (for instance,  condition (\ref{cov_max})
is the same as condition (2.40) of Arcones (1994) in the case of stationary sequences).
We expect that parts (i) and (ii) can be also obtained following the method in Nourdin {\it et al.} (2010).
Part (iii) seems more interesting. In Bardet and Surgailis (2011), (iii) is applied when ${\mbf X}_{\! n}(j)= {\mbf Z}_{j/n}$ and $({\mbf Z}_t)_t$ is a vector valued continuous time process. \\

Similarly to Lemma \ref{lemgauss}, Theorem \ref{tlcgauss} can be extended to nonstandardized Gaussian vectors.
Corollary \ref{tlcgauss2} below refers to the most interesting part (iii) of  Theorem \ref{tlcgauss}.

%A particular case ($m=2$) for non-standardized Gaussian vectors is proposed below:

\begin{cor}\label{tlcgauss2}
Let ${\mbf
Y}_{\! n}(k) =  \big (Y^{(1)}_{\! n}(k), \dots, Y^{(\nu)}_{\! n}(k)\big ) \in \R^\nu, \, 1 \le k \le n, \, n \in \N $ be a  triangular
array of jointly Gaussian vectors,  with zero mean $\E {\mbf Y}_{\! n}(k)= 0$ and non-degenerate covariance matrices $\Sigma_{k,n} =
\E  {\mbf
Y}_{\! n}(k) {\mbf
Y}_{\! n}(k)^\intercal.$ 
Assume that covariances 
$r^{(p,q)}_n(j,k) :=  {\rm Cov}\big(Y^{(p)}_n(j), Y^{(q)}_n(k)\big) $ satisfy (\ref{rho_dom}), for some $m \ge 1$.
Moreover, assume that (\ref{Y_W}) holds with ${\mbf
X}_{\! n}(\cdot)$ replaced by ${\mbf
Y}_{\! n}(\cdot),$ where $ \left({\mbf W}_{\! \tau}(j)\right)_{j\in {\Z}} $ is 
a stationary Gaussian ${\R}^\nu$-valued process with non-degenerate covariance matrix
$\Sigma_\tau := \E {\mbf W}_{\! \tau}(0) {\mbf W}_{\! \tau}(0)^\intercal $
such that
\begin{equation} \label{bddS}
\sup_{\tau \in (0,1]} \| \Sigma^{-1}_\tau \| < \infty
\end{equation}
and 
\begin{equation} \label{approxS}
\sup_{\tau \in (0,1]}\|\Sigma_{[n \tau],n} - \Sigma_\tau \| \limiten   0.
\end{equation}
Let $f_{k,n} \in \L^2_0({\mbf
Y}_{\! n}(k)), \, 1\le k \le n, \, n \in \N$ be a triangular array of functions all having a generalized  Hermite rank
not less than $m$ and such that
\begin{eqnarray} \label{rank22}
\sup_{\tau \in (0,1]} \E \left (\tilde f_{[\tau n],n}({\mbf X}) - \tilde \phi_{\tau}({\mbf X})\right )^2 \limiten 0,
\end{eqnarray}
where $\tilde f_{k,n}({\mbf x}) :=f_{k,n}( \Sigma^{1/2}_{k,n} {\mbf x}) $ and where
$ \tilde \phi_\tau, \tau \in [0,1]$ is a $\L^2_0({\mbf X})-$valued continuous function, with ${\mbf X}$ a standard Gaussian vector
in $\R^\nu$ as usual.
Then 
\begin{equation} \label{CLTbis}
n^{-1/2} \sum_{k=1}^n  f_{k,n} \big({\mbf Y}_{\! n}(k) \big) % - \E  f_{k,n} \big({\mbf Y}_{\! n}(k)\big)   \Big )
\limiteloin {\mathcal N}(0, \sigma^2).
\end{equation}
where $ \sigma^2$ is defined  in (\ref{sigma}), with  $\phi_\tau ({\mbf x}) := \tilde \phi_\tau (\Sigma^{-1/2}_\tau {\mbf x})$.

\end{cor}

\noindent {\it Proof of Corollary \ref{tlcgauss2}.} Similarly as in the proof of Corollary \ref{corgauss},
let ${\mbf X}_{\! n}(k) :=  \Sigma^{-1/2}_{k,n} {\mbf Y}_{\! n}(k)$. % \,  \tilde f_{k,n}({\mbf x}) := f_{k,n}( \Sigma^{1/2}_{k,n} {\mbf x}).$
The  ${\mbf X}_{\! n}(k)$'s are standardized Gaussian vectors in $\R^\nu$ and the $\tilde f_{k,n}$'s belong to $\L^2({\mbf X})$ and have
a Hermite rank not less than $m$.
%are measurable functions with the property that
%$\E [ {\mbf Y}_{\! n}(k)  \tilde f_{k,n}({\mbf Y}_{\! n}(k))] =  0 $
%implying that the Hermite rank of $\big (\tilde f_{k,n}- \E \tilde f_{k,n}({\mbf Y}_{\! n}(k))\big )$ is not less than $2$.
Assumptions
(\ref{Y_W}) and (\ref{approxS}) entail for any $\tau \in (0,1)$, $\big({\mbf X}_{\! n}(j + [n\tau]))
\big)_{-J \le j \le J}$
$\limiteloin $  $ \big(\widetilde {\mbf W}_{\! \tau}(j)\big)_{-J \le j \le J}$, where
$\widetilde {\mbf W}_{\! \tau}(j)
 := \Sigma^{-1/2}_\tau {\mbf W}_{\! \tau}(j),  \, j\in {\Z}$ is a
stationary Gaussian process having a unit covariance matrix $\E \widetilde {\mbf W}_{\! \tau}(0) \widetilde {\mbf W}_{\! \tau}(0)^\intercal = I$.
Conditions (\ref{bddS}) and (\ref{approxS})  imply that $\max_{1\le k \le n} \|\Sigma_{k,n}^{-1/2}\| \le C$.
The last fact together with condition (\ref{rho_dom}) for covariances $r^{(p,q)}_n(j,k) :=  {\rm Cov}\big(Y^{(p)}_n(j), Y^{(q)}_n(k)\big) $
imply a similar condition for ${\rm Cov}\big(X^{(p)}_n(j), X^{(q)}_n(k)\big)$:
for all $(j, k)\in \{1,\cdots,n\}^2$ we have that
$\max_{1\le u, v \le \nu} \big|
\E X^{(u)}_{\! n}(j) X^{(v)}_{\! n}(k) \big| \le C|\rho(j-k)|$; see (\ref{XYbdd}).
% Moreover, with ${\mbf X}$ a standardized Gaussian vector and  $\tilde \phi_\tau({\mbf x})=\phi_\tau (\Sigma^{1/2}{\mbf x})$, condition %(\ref{approxS})
%ensures that $\sup_{0\leq \tau\leq 1 } n \, \E\big ( \tilde f_{[n\tau],n}({\mbf X} )- \tilde \phi_\tau({\mbf X})\big )^2  \limiten 0 $.
This way we see that the conditions of Theorem \ref{tlcgauss}(iii) including (\ref{f_n})
are satisfied and can be applied to the families of Gaussian vectors
$({\mbf X}_{\! n}(k))$ and functions $\big (\tilde f_{k,n}\big )$, yielding (\ref{CLTbis}).  \hfill $\Box$ \\

\medskip

\noindent {\it Proof of Theorem \ref{tlcgauss}.} (i) Using Arcones' inequality (see Arcones, 1994, (2.44) or Soulier, 2001, (2.4)), one obtains
\begin{eqnarray*}
\E \Big(n^{-1/2} \sum_{k=1}^n f_k\left({\mbf X}_{\! n}(k)\right)
\Big)^2 & =& \frac 1 n \sum_{k=1}^n \|f_k\|^2 + \frac 1 n
\sum\nolimits^\prime \E f_k({\mbf X}_{\! n}(k))f_\ell({\mbf X}_{\!
n}(\ell)) \\
& \leq & \max_{1\leq k \leq n} \|f_k\|^2 +C \big (\max_{1\leq k
\leq n} \|f_k\|\big )^2 \max_{1\leq k \leq n} \sum_{1\le \ell \le
n, \ell\neq k} \max_{1\le p, q\le \nu} \big
|r^{(p,q)}_n(k,\ell)\big |^m,
\end{eqnarray*}
%thanks to Lemma \ref{lemgauss} (with $\alpha=p=2$),
where $C$ is a
positive real number not depending on $n$ or $f_k$. Now, using
assumption (\ref{cov_max}), (i) is proved.
\medskip

\noindent (ii)  We use the following well-known fact. Let $(Z_n)_{n
\ge 1} $ be a sequence of r.v.'s with zero mean and finite variance.
Then $Z_n   \limiteloin  {\mathcal N}(0, \sigma^2)$ if and only if
for any $\epsilon>0$ one can find an integer $n_0(\epsilon) \ge 1$
and a sequence $(Z_{n,\epsilon})_{n\ge 1}$ satisfying
$Z_{n,\epsilon} \limiteloin {\mathcal N}(0, \sigma^2_\epsilon)$ and
$\forall n > n_0(\epsilon)$, $ \E (Z_n - Z_{n,\epsilon})^2 <
\epsilon$.

 Let $Z_n := n^{-1/2} \sum_{k=1}^n f_{k,n}\left({\mbf
X}_{\! n}(k)\right) $. We shall construct an approximating sequence
$Z_{n,\epsilon}$ with the above properties in two steps.

Firstly, by condition (\ref{f_n}) and continuity of $\phi_\tau $,
for a given $\epsilon> 0$ one can find integers $M, n_0(\epsilon)
$ and a partition $0=:\tau_0 < \tau_1 < \dots < \tau_M <
\tau_{M+1} := 1$ such that $\forall \  n > n_0(\epsilon)$,
\begin{equation}
\max_{0\le i \le M} \max_{k/n \in (\tau_i, \tau_{i+1}]}\| f_{k,n}
- \phi_{\tau_i}\| =\max_{0\le i \le M} \max_{k/n \in (\tau_i,
\tau_{i+1}]}\big (\E ( f_{k,n}({\mbf X}) - \phi_{\tau_i}({\mbf
X}))^2\big) ^{1/2} \ < \ \epsilon. \label{norm}
\end{equation}
Put
$$
\widetilde Z_{n,\epsilon} := n^{-1/2} \sum_{i=0}^M \sum_{k/n \in
(\tau_i, \tau_{i+1}]} \phi_{\tau_i}\left({\mbf X}_{\!
n}(k)\right).
$$
Note for any $\tau \in (0,1]$, the function $\psi_\tau $ has
Hermite rank not less than $m$, being the limit of a sequence of $
\L^2_0({\mbf X})-$valued functions of Hermite rank $\ge m$.
Therefore for the difference $Z_n - \widetilde Z_{n,\epsilon}$ the
inequality (\ref{f_bdd}) applies, yielding $\forall \  n>
n_0(\epsilon)$
\begin{equation}
\E (Z_n - \widetilde Z_{n,\epsilon})^2 \ \le \  C \max_{0\le i \le
M} \max_{k/n \in (\tau_i, \tau_{i+1}]}\| f_{k,n} -
\phi_{\tau_i}\|^2 \leq   C \epsilon^2 \label{Z1approx}
\end{equation}
in view of (\ref{norm}), with a constant $C$ independent of $n,
\epsilon $.

Secondly, we expand each $\phi_{\tau_i}$ in Hermite polynomials:
\begin{eqnarray}
\phi_{\tau_i}({\mbf x})&=& \sum_{m \le |{\mbf k}|} \frac{J_i({\mbf
k})}{{\mbf k}!} H_{\mbf k}({\mbf x}),  \qquad (i=0,1, \dots, M)
\end{eqnarray}
where
$$
J_i({\mbf k}) := J_{\phi_{\tau_i}}({\mbf k}) =  \E
\phi_{\tau_i}({\mbf X}) H_{\mbf k}({\mbf X}), \qquad |J_i({\mbf k})|
\le   \|\phi_{\tau_i}\| ({\mbf k}!)^{1/2}.
$$
We can choose $t(\epsilon) \in {\N}$ large enough so that
\begin{equation}
\|\phi_{\tau_i} - \phi_{\tau_i, \epsilon}\| \le \epsilon, \qquad
(i=0,1, \dots, M), \label{phi_approx}
\end{equation}
where $\phi_{\tau_i, \epsilon}$ is a finite sum of Hermite
polynomials:
\begin{equation}
\phi_{\tau_i, \epsilon}({\mbf x}) \ := \  \sum_{m \le |{\mbf k}|
\le t(\epsilon)} \frac{J_i({\mbf k})}{{\mbf k}!} H_{\mbf k}({\mbf
x}),  \qquad (i=0,1, \dots, M).
\end{equation}
Note $t(\epsilon) $ does not depend on $i=0,1,\dots, M$, and
$\epsilon>0$ is the same as in (\ref{norm}). Put
\begin{equation}
Z_{n,\epsilon} := n^{-1/2} \sum_{i=0}^M \sum_{k/n \in (\tau_i,
\tau_{i+1}]} \phi_{\tau_i, \epsilon}\left({\mbf X}_{\!
n}(k)\right).  \label{Z_epsilon}
\end{equation}
Applying (\ref{f_bdd}) to the difference $\widetilde
Z_{n,\epsilon} - Z_{n,\epsilon}$ and using (\ref{phi_approx}) and
(\ref{Z1approx}), we obtain $\forall \  n> n_0(\epsilon)$,
\begin{equation}
\E (Z_n - Z_{n,\epsilon})^2 \  \le \   C \epsilon^2
\label{Z2approx}
\end{equation}
where the constant $C$ is independent of $n, \epsilon $.  Let
$\sigma^2_{n,\epsilon} :=  \E Z^2_{n,\epsilon}$. From
(\ref{Z2approx}) and condition (\ref{sigma_n}) it follows that
$\forall n>n_0(\epsilon)$,
\begin{equation}
\sigma^2 - C \epsilon \  \le \  \sigma^2_{n,\epsilon}  \ \le \
\sigma^2 + C \epsilon,
 \label{sigma_sigma}
\end{equation}
with some $C$ independent of $n, \epsilon $. In particular, by
choosing $\epsilon >0$ small enough, it follows that $\liminf_{n
\to \infty} \sigma^2_{n,\epsilon} > 0$. We shall prove below that
for any fixed $\epsilon >0$,
\begin{equation}
U_n \ := \  \frac{Z_{n,\epsilon} }{\sigma_{n,\epsilon}} \ =  \
\frac{1} {\sigma_{n, \epsilon} n^{1/2}} \sum_{i=1}^M \sum_{k/n \in
(\tau_i, \tau_{i+1}]} \phi_{\tau_i, \epsilon} \left({\mbf X}_{\!
n}(k)\right) \limiteloin {\mathcal N}(0,1).
 \label{CLT2}
\end{equation}
As noted in the beginning of the proof of the theorem, the CLT in
(\ref{CLT}) follows from (\ref{CLT2}), (\ref{Z2approx}),
(\ref{sigma_sigma}). Indeed, write
\begin{eqnarray*}
\E \e^{\i a Z_n} - \e^{-a^2 \sigma^2/2} &=&\left(\E \e^{\i a Z_n}-
\E \e^{\i a Z_{n,\epsilon}}\right) +
\left(\E \e^{\i a \sigma_{n,\epsilon} U_{n}} - \e^{-a^2\sigma^2_{n,\epsilon}/2}\right) \\
&+&\left( \e^{-a^2\sigma^2_{n,\epsilon}/2} -
\e^{-a^2\sigma^2/2}\right) \ : = \ \sum_{i=1}^3 \ell_i(n).
\end{eqnarray*}
Here, for some constant $C$ independent of $n, a, \epsilon $,
\begin{eqnarray*}
|\ell_1(n)|&\le&\E^{1/2}\big|{\e}^{\i a(Z_n-Z_{n,\epsilon})}
-1\big |^2 \
\le \ |a|\E^{1/2} |Z_n -Z_{n,\epsilon}|^2 \ \le\  C|a| \epsilon, \\
|\ell_3(n)|&\le&Ca^2 \left|\sigma^2_{n,\epsilon} - \sigma^2\right|
\ \le \  Ca^2 \epsilon,
\end{eqnarray*}
and therefore $\ell_i(n), i =1, 3 $ can be made arbitrarily small by
choosing $\epsilon >0$ small enough; see (\ref{Z2approx}),
(\ref{sigma_sigma}).  On the other hand, the convergence in
(\ref{CLT2}) implies  uniform convergence of characteristic
functions on compact intervals and therefore $\sup_{|a| \le A}
|\ell_2(n)| \le \sup_{|a|\le 2A}  \left|\E \e^{\i a U_{n}} -
\e^{-a^2/2}\right| \limiten 0$ for any $A>0$.  This proves
(\ref{CLT}).

\smallskip

 It remains to prove (\ref{CLT2}). The proof of the
corresponding CLTs for sums of Hermite polynomials in Arcones (1994)
and Breuer and Major (1983) refer to stationary processes and use
Fourier methods. Therefore we present an independent proof of
(\ref{CLT2}) based on cumulants and the H\"older inequality in
(\ref{H_bdd}). Again, our proof appears to be much simpler than
computations in the above mentioned papers.

Accordingly, it suffices to show that cumulants of order $p \ge 3$
of \ $U_n$ \  asymptotically vanish. In view of
(\ref{sigma_sigma}) and linearity of cumulants, this follows from
the fact that for any $p\ge 3$ and any multiindices ${\mbf k}_u =
(k^{(1)}_u, \dots, k^{(\nu)}_u) \in {\Z}^\nu_+, \  u=1, \dots, p$
with $k_u = |{\mbf k}_u| = k^{(1)}_u+ \dots + k^{(\nu)}_u \ge m \
(1\le u \le p)$,
\begin{equation}
\Sigma_n \ := \ \sum_{t_1, \dots, t_p =1}^n |{\rm cum}(t_1, \dots,
t_p)|   = \  o(n^{p/2}),  \label{cum_o}
\end{equation}
where ${\rm cum}(t_1, \dots, t_p)$ stands for joint cumulant:
\begin{equation}
{\rm cum}(t_1, \dots, t_p) \ := \ {\rm cum}\left(H_{{\mbf
k}_1}({\mbf X}_{\! n}(t_1)), \dots,  H_{{\mbf k}_p}({\mbf X}_{\!
n}(t_p))   \right). \label{cum}
\end{equation}
Split  $\Sigma_n = \Sigma'_n(K) + \Sigma''_n(K), $ where
$$
\Sigma'_n(K)  \  := \ \sum_{t_1, \dots,t_p=1}^n  \left|{\rm
cum}(t_1,\dots,t_p)\right| {\mbf 1}(|t_i-t_j| \le K \  \forall
i\ne j)
$$
and where $K $ will be chosen large enough. Then for any fixed
$K$, we have $\Sigma'_n(K) = O(n) = o(n^{p/2}) $ as $p\ge 3$. The
remaining sum $\Sigma''_n(K)$ does not exceed $\sum_{1\le i \ne j
\le p} \Sigma''_{n,i,j}(K)$, where
$$
\Sigma''_{n,i,j}(K)\ := \ \sum_{t_1, \dots, t_p=1}^n \left|{\rm
cum}(t_1,\dots,t_p)\right| {\mbf 1}(|t_i-t_j| > K).
$$
Therefore, relation (\ref{cum_o}) follows if we show that there
exist $\delta(K)  \limiteK 0$ and $\tilde n_0$ such that for any
$1\le i\ne j \le p$ and any $n> \tilde n_0$
\begin{equation}
\limsup_{n\to \infty} \Sigma''_{n,i,j}(K) \ < \  \delta(K) n^{p/2}.
\label{delta_K}
\end{equation}
The proof below is limited to $(i,j)=(1,2)$ as the general case is
analogous. It is well-known that the joint cumulant in
(\ref{cum}), similarly to the joint moment in (\ref{hermite_bdd}),
can be expressed as a sum over all {\em connected} diagrams
$\gamma \in \Gamma_c(T)$ over the table $T$ in (\ref{T_table}). By
introducing $\bar \rho (s,t) := \max_{1\le p,q \le \nu}
\left|r^{(p,q)}_n(s,t)\right|$,  we obtain
\begin{equation}
|{\rm cum}(t_1, \dots, t_p)| \ \le \ \sum_{\gamma \in \Gamma_c(T)}
\prod_{1 \le u <v \le p} (\bar \rho (t_u,t_v))^{\ell_{uv}},
\label{T1}
\end{equation}
where we use the notation in (\ref{hermite_bdd}). Therefore,
$$
\Sigma''_{n,1,2}(K) \  \le \  \sum_{\gamma \in \Gamma_c(T)}
\sum_{t_1, \dots, t_p =1}^n \prod_{1 \le u <v \le p} (\bar \rho
(t_u,t_v))^{\ell_{uv}} {\mbf 1}(|t_1-t_2|>K) \ := \ \sum_{\gamma
\in \Gamma_c(T)} \bar I_{n,T}(\gamma),
$$
Next, by applying the H\"older inequality as in (\ref{H_bdd}),
\begin{eqnarray}
\bar I_{n,T}(\gamma) &\le& \min\left(\prod_{1\le u<v\le p} \bar
R_{uv}, \prod_{1\le u<v\le p} \bar R_{vu} \right). \label{H_bdd1}
\end{eqnarray}
where (cf. (\ref{R_uv}))
\begin{eqnarray*}
\bar R_{uv} &:=&\cases{\left(\sum\nolimits_{1\le t\le n}
\left(\sum\nolimits_{1\le s \le n} \bar \rho^{k_u}(s,t)
\right)^{k_v/k_u} \right)^{\ell_{uv}/k_v}, &$(u,v) \ne (1,2),
(2,1)$, \cr \left(\sum\nolimits_{1\le t\le n}
\left(\sum\nolimits_{1\le s \le n} \bar \rho^{k_1}(s,t) {\mbf
1}(|t-s|>K) \right)^{k_2/k_1}
\right)^{\ell_{12}/k_2},&$(u,v)=(1,2)$,\cr
\left(\sum\nolimits_{1\le t\le n} \left(\sum\nolimits_{1\le s \le
n} \bar \rho^{k_2}(t,s) {\mbf 1}(|t-s|>K) \right)^{k_1/k_2}
\right)^{\ell_{12}/k_1},&$(u,v)=(2,1).$\cr}
\end{eqnarray*} From assumptions  (\ref{cov_max}), (\ref{cov_rank_K}), there exists
a constant $C $ and $\delta(K) \limiteK 0$ independent of $n$ such
that for any $k\ge m$ and any $n\ge 1$
\begin{eqnarray*}
\sup_{1\le t \le n} \sum_{s=1}^n \bar \rho^k (s,t) &\le& Cn, \\
\sup_{1\le t \le n} \sum_{s=1}^n \bar \rho^k (s,t){\mbf
1}(|t-s|>K)&\le&\delta(K)n.
\end{eqnarray*}
Therefore
$$
\bar R_{uv} \  \le  \ \cases{Cn^{\ell_{uv}/k_v},  &$(u,v) \ne
(1,2),  (2,1)$, \cr \tilde
\delta(K)n^{\ell_{12}/k_2},&$(u,v)=(1,2)$,\cr \tilde
\delta(K)n^{\ell_{12}/k_1},&$(u,v)=(2,1)$,\cr }
$$
with some $\tilde \delta(K) \limiteK 0 $ independent of $n$.
Consequently, the minimum on the right-hand side of (\ref{H_bdd1})
does not exceed
$$
C \tilde \delta(K) \min\left(n^{\sum_{1\le u < v\le p}
\ell_{uv}/k_v}, n^{\sum_{1\le u < v\le p} \ell_{uv}/k_u}\right) \
=  \  C \tilde \delta(K) n^{\min (L(T), L^*(T))}
$$
where the quantities $L(T), L^*(T)$ introduced in (\ref{L(T)})
satisfy $L(T) + L^*(T) = p$, see (\ref{U_iden}), and therefore
$\min(L(T), L^*(T)) \le p/2$. This proves  (\ref{delta_K}) and the
CLT in (\ref{CLT2}), thereby   completing the
proof of part (ii).

\smallskip

\noindent (iii) Let us first prove (\ref{sigma_n}) with $\sigma^2
$ given in (\ref{sigma}) in the case when $f_{k,n} \equiv f$ do
not depend on $k, n$ (in such case, one has $\phi_\tau \equiv f $,
too). We have
\begin{eqnarray}
\sigma^2_n &=&n^{-1}\sum_{k,k'=1}^n \E \big[f\left({\mbf X}_{\!
n}(k)\right) f\left({\mbf X}_{\! n}(k')\right)\big] \ =  \  \int_0^1
F_n(\tau) \d \tau, \label{int_F}
\end{eqnarray}
where
\begin{equation}
F_n(\tau) \ := \  \sum_{j=1- [n\tau]}^{n - [n\tau]} \E \big[
f\left({\mbf X}_{\! n}([n\tau])\right) f\left({\mbf X}_{\!
n}([n\tau] + j)\right)\big].  \label{F_n}
\end{equation}
Condition (\ref{Y_W}) implies that
$$
\E \big[f\left({\mbf X}_{\! n}([n\tau])\right) f\left({\mbf X}_{\!
n}([n\tau] + j)\right)\big] \rightarrow \E \big[f\left({\mbf W}_{\!
\tau}(0)\right) f\left({\mbf W}_{\! \tau}(j)\right)\big] $$ for each
$j\in {\Z}$ as $n \to \infty $. From (\ref{rho_dom}) and with the
inequality of previous part (i), it exists $C>0$ such that
\begin{equation}
\Big|\E \big[ f\left({\mbf X}_{\! n}([n\tau])\right) f\left({\mbf X}_{\!
n}([n\tau] + j)\right) \big] \Big |\ \leq \  C |\rho(j)|^m, \label{arcones}
\end{equation}
and $\sum_{j\in {\Z}} |\rho(j)|^m < \infty $. Hence, from Lebesgue
Theorem,
$$
F_n(\tau)=\sum_{j \in \Z} \1_{j \in \{ 1- [n\tau],\cdots,n -
[n\tau]\}}\E \big[f\left({\mbf X}_{\! n}([n\tau])\right) f\left({\mbf
X}_{\! n}([n\tau] + j)\right)\big]  \limiten  \sum_{j \in \Z}\E
\big[f\left({\mbf W}_{\! \tau}(0)\right) f\left({\mbf W}_{\!
\tau}(j)\right)\big].
$$
The dominated convergence theorem allows one to  pass to the limit
under the integral, thereby proving (\ref{sigma_n}) with $\sigma^2$
given in (\ref{sigma}) in the case $f_{k,n} \equiv f$.

\smallskip

To end the proof, consider the general case of $f_{k,n}$ as in
(iii). Let $Z_{n,\epsilon}$ be defined as in (\ref{Z_epsilon}). Note
relation (\ref{Z2approx}) holds as its proof does not use
(\ref{sigma_n}). In part (ii), we used  (\ref{sigma_n}) to prove
(\ref{sigma_sigma}). Now we want to prove  (\ref{sigma_sigma}) using
(\ref{Y_W}) instead of (\ref{sigma_n}). This will suffice for the
proof of (iii), as the remaining argument is the same as in part
(ii).

Consider the variance  $\sigma^2_{n,\epsilon} = \E
Z^2_{n,\epsilon}$ of $Z_{n,\epsilon} $ defined  in
(\ref{Z_epsilon}):
\begin{eqnarray*}
\sigma^2_{n,\epsilon} &=&n^{-1}\left(\sum_{0\le i \le M} \E D^2_i
+ 2\sum_{0\le i < j \le M} \E D_i D_j\right),
\end{eqnarray*}
where
$$
D_i \  := \   \sum_{k/n \in [\tau_i, \tau_{i+1})} \phi_{\tau_i,
\epsilon}\left({\mbf X}_{\! n}(k)\right).
$$
Let us show that for $\epsilon, M $ fixed, and as $n \to \infty $,
\begin{eqnarray}
\E D_i D_j &=&o(n)\qquad \qquad (i \neq j),  \label{D_ij}   \\
n^{-1}\E D_i^2&\limiten &\int_{\tau_i}^{\tau_{i+1}} \sum_{j\in
{\Z}} \E \big[\phi_{\tau_i,\epsilon} \left({\mbf W}_{\! \tau}(0)\right)
\phi_{\tau_i,\epsilon} \left({\mbf W}_{\! \tau}(j)\right) \big]\d \tau.
\label{D_i}
\end{eqnarray}
Here, (\ref{D_i}) follows from the argument in the beginning of the
proof of (iii),  as $\phi_{\tau_i,\epsilon}$ does not depend on $k,
n$.  Relation (\ref{D_ij}) is implied by the following computations.
Using the Hermitian rank of functions $\phi_{\tau_i, \epsilon}$, for
$i <j$ one obtains
\begin{eqnarray*}
\left|\E  \phi_{\tau_i, \epsilon}\left({\mbf X}_{\!
n}([n\tau_i]+k)\right) \phi_{\tau_j, \epsilon}\left({\mbf X}_{\!
n}([n\tau_j]+\ell)\right) \right| &\le & C \| \phi_{\tau_i,
\epsilon} \| \cdot  \| \phi_{\tau_j, \epsilon} \|\max_{1\le p,
q\le \nu} \left|r^{(p,q)}_n([n\tau_i]+k,[n\tau_j]+\ell)\right|^m\\
& \leq & C \| \phi_{\tau_i, \epsilon} \| \cdot  \| \phi_{\tau_j,
\epsilon} \|  \big | \rho ([n\tau_j]-[n\tau_i]+\ell-k)\big |^m.
\end{eqnarray*}
Therefore, for $i<j$, and $\varepsilon$ small enough,
\begin{eqnarray*}
\big | \E D_i D_j \big | & \leq & C \max_{\tau\in [0,1]} \|
\phi_{\tau} \|^2 \sum_{k=0}^{[\tau_{i+1}n]-[\tau_in]}
\sum_{\ell=0}^{[\tau_{j+1}n]-[\tau_jn]} \big | \rho
([n\tau_j]-[n\tau_i]+\ell-k)\big |^m \\
& \leq & C \max_{\tau\in [0,1]} \| \phi_{\tau} \|^2 \sum_{k=1}^n k
\big | \rho(k)\big |^m \ = \  o(n)
\end{eqnarray*}
since $\sum_{k=1}^n k \big | \rho(k)\big |^m \le \sqrt{n} \sum_{1\le
k \le \sqrt{n}} \big | \rho(k)\big |^m + n \sum_{k > \sqrt{n}} \big
| \rho(k)\big |^m = o(n)$. 
Thus, (\ref{D_ij}) is proved. From (\ref{D_ij}), (\ref{D_i}) it
follows that for any $\epsilon>0$
$$
\lim_{n \to \infty} \sigma^2_{n, \epsilon} \ =  \  \bar
\sigma^2_\epsilon \ := \ \sum_{i=0}^M \int_{\tau_i}^{\tau_{i+1}}
\sum_{j\in {\Z}} \E \big[\phi_{\tau_i,\epsilon} \left({\mbf W}_{\!
\tau}(0)\right) \phi_{\tau_i,\epsilon} \left({\mbf W}_{\!
\tau}(j)\right)\big] \d \tau.
$$
Consider the difference \  $ \bar \sigma^2_\epsilon - \sigma^2 =
\sum_{i=0}^M \int_{\tau_i}^{\tau_{i+1}} \sum_{j\in {\Z}}
\Theta_{M,\epsilon}(\tau,j) \d \tau, $ \  where
\begin{eqnarray}
|\Theta_{M,\epsilon}(\tau,j)|&=&\left|\E  \phi_{\tau_i,\epsilon}
({\mbf W}_{\! \tau}(0)) \phi_{\tau_i,\epsilon}({\mbf W}_{\!
\tau}(j))
- \E \phi_{\tau} ({\mbf W}_{\! \tau}(0)) \phi_{\tau}({\mbf W}_{\! \tau}(j)) \right| \nonumber  \\
&\le&\left|\E \left(\phi_{\tau_i,\epsilon}({\mbf W}_{\! \tau}(0))
- \phi_{\tau}({\mbf W}_{\! \tau}(0)) \right)
\phi_{\tau_i,\epsilon}({\mbf W}_{\! \tau}(j)) \right| \ +\
\left|\E \left(\phi_{\tau_i,\epsilon} \left({\mbf W}_{\!
\tau}(j)\right) - \phi_{\tau} \left({\mbf W}_{\!
\tau}(j)\right)\right)
\phi_{\tau} \left({\mbf W}_{\! \tau}(0) \right)\right| \nonumber \\
&\le&\|\phi_{\tau_i,\epsilon} - \phi_\tau\| \left(\|\phi_{\tau_i,
\epsilon}\| + \|\phi_{\tau}\| \right). \label{U_bdd}
\end{eqnarray}
Using uniform continuity of $\phi_\tau, \tau \in [0,1]$ (in the
sense of $\L^2-$norm continuity), we obtain that the right-hand side
of (\ref{U_bdd}) can be made arbitrarily small by choosing $M $ (=
the number of partition intervals of $[0,1]$) and $t(\epsilon)$ (=
the truncation level of Hermite expansion) sufficiently large,
uniformly in $\tau \in [0,1]$ and $j\in {\Z}$. On the other hand, $
|\Theta_{M,\epsilon}(\tau,j)| \le C\sup_{\tau \in [0,1]}\|\phi_\tau
\|^2 |\rho(j)|^m $
 by Arcones' inequality, c.f. (\ref{arcones}). Therefore
 $|\Theta_{M,\epsilon}(\tau,j)|$ is dominated by a summable function
 uniformly in $M, \epsilon$. Now, (\ref{sigma_sigma}) follows by an application of Lebesgue
 theorem. This proves part (iii) and Theorem \ref{tlcgauss} too. \hfill $\Box $

 \section{A Berry-Esseen-type bound for nonstationary Gaussian subordinated triangular arrays}\label{breuer}

This section obtains a Berry-Esseen-type upper bound in the
CLT (\ref{CLT}) for non-stationary Gaussian subordinated triangular arrays
following %can be specified using an extension of
the method and results presented in
Nourdin {\it et al.} (2011).
We will refer NPP to the last paper  in the rest of this section.
To simplify the discussion,
%Indeed this paper is devoted to the stationary case and we are going to extend it to the case our non-stationary case of triangular arrays.
we restrict our task to the case when the functions $f_{k,n} = f$  in
Theorem \ref{tlcgauss} (iii) do not depend on $k, n$.
%(the case of a family $(f_{k,n})$ requires too much changes w.r.t. the results in Nourdin {\it et al.}: for instance the crucial relationship (4.41)  %does not hold anymore).
As in NPP, our starting point is the Hermite expansion (\ref{hermite0}) written as
\begin{equation} \label{hermiteBE}
f = \sum_{\ell = m}^\infty f_{(\ell)}, \qquad f_{(\ell) } :=  \sum_{|{\mbf k}|= \ell} J_f({\mbf k}) H_{\mbf k}/{\mbf k}!.
\end{equation}
Following NPP and using the Hermite expansion in (\ref{hermiteBE}),
%of $f=\sum_{\ell=m}^\infty f_\ell$ (see (2.7) in Nourdin {\it et al.}),
we first define  the following quantities: %already defined in Nourdin {\it et al.}
for $j\in \Z$, $\ell\geq m$, $N\geq m$, $n\in \N^*$ and $J\in \{1,\dots,n\}$:
\begin{eqnarray}
\theta(j)&:=& |\rho(j)|~,\quad
K:=\inf\{k\in \N: \theta(j) \leq \frac 1 \nu ,\, \forall |j|\geq k \big \},  \quad
\theta :=  \sum_{j\in \Z} \theta(j)^m, \\
\label{sigmam}\sigma_{\ell,n}^2&:=&n^{-1}\sum_{t, t'=-n}^n  \Cov\big (f_{(\ell)}\left({\mbf
X}_{\! n}(t) \right), f_{(\ell)}\left({\mbf
X}_{\! n}(t') \right) \big),  \\
\label{gamma}  \gamma_{n,\ell,e}&:=& \frac 1 {n^{1/2}} \, \Big (2\theta \,
 \sum_{|j|\leq n} \theta(j)^e \sum_{|j'|\leq n} \theta(j')^{\ell-e}\Big )^{1/2} \,\,\, (\mbox{for } 1\leq e\leq \ell-1), \\
%\label{A1}A_{1,n}&=& \frac {\E[f^2({\mbf X})] }{2} \Big [ \frac{2K^2}{n}+ \nu^q\big(\sum_{|j|\leq n} \theta(j)^m \frac{|j|}{n} +\sum_{|j|>n}\theta(j)^m\big) \Big]\\
\label{A2} A_{2,N}&:=& 2(2K+\nu^m \theta)
\Big (\E[f^2({\mbf X})]\sum_{\ell=N+1}^\infty \E[f^2_{(\ell)}({\mbf X})]\Big)^{1/2} ,\\
% \label{A3} A_{3,n,N}\!\!&=&\!\! \frac{\E[f^2({\mbf X})]}{2}\sum_{\ell=m}^N \Big (\frac{\nu^\ell }{\ell \,  \ell!} \sum_{k=1}^{\ell-1}k \, k!{\ell \choose k}^2\sqrt{(2\ell-2k )!} \,\gamma_{n,\ell,k}\Big) \\
\label{A3}A_{3,n,N}&:=&\frac{1}{2} \E [f^2 ({\mbf X})] \sum_{\ell = m}^N
\Big( \frac{\nu^\ell}{ \ell \ell !} \sum_{j=1}^{\ell -1} j j! {\ell \choose j}^2 \sqrt{(2\ell - 2j)!} \gamma_{n,\ell, j} \Big), \\
\label{A4} A_{4,n,N}&:=&\frac1 2 \, \E[f^2({\mbf X})] %(2K+\nu^m \theta)^{1/2} \!
\sum_{m\leq \ell <\ell'\leq N} \nu ^{{\ell'}/2}\sqrt{\frac{\ell'!}{\ell!}}
\frac{\ell +\ell'}{\ell}{\ell'-1 \choose \ell-1}\big ((\ell'-\ell)! \gamma_{n,\ell',\ell'-\ell}\big ) ^{1/2}, \\
 A_{5,n,N}&:=& \frac{\E[f^2({\mbf X})]}{2\sqrt{2}}\!\!\!\! \!\!\!\!\!  \sum_{m\leq \ell<\ell'\leq N}\!\!\!\!\!\!\!  (\ell +\ell')
\sum_{j=1}^{\ell-1} \!(j-1)!{\ell\!\!-\!\!1 \choose j\!\!-\!\! 1}{\ell'\!\!-\!\!1 \choose j\!\!-\!\!1}\sqrt{(\ell\!+\!\ell'\!-\!2j)!} \Big (\frac{\nu^\ell}{\ell!}
\gamma_{n,\ell,\ell-j}\!\!+\!\!\frac{\nu^{\ell'}}{\ell'!} \gamma_{n,\ell',\ell'-j}\Big), \label{A5}\\
A_{6,n,J}&:=&\frac 1 2 |\partial f|^2_\infty \sup_{0\le \tau \le 1} \sum_{|j| \le J} \Big\| \E [{\mbf
X}_{\! n}([n\tau]) {\mbf
X}^\intercal_{\! n}([n\tau] +j) ] - \E [ {\mbf W}_\tau(0) {\mbf W}^\intercal_\tau(j) ]\Big\|, \label{A6}\\
%\Big ( \E \big [ f^2({\mbf X})\big ] \,\| \partial f \|_\infty^2 \Big )^{1/2} \, (2N+1)\, \sqrt{ u_N(n)}
A_{7,J}&:=&\frac 1 2  \E \big [ f^2({\mbf X})\big ]  \nu^m  \sum_{|k|>J} \theta^m(k)   \label{A7}.
\end{eqnarray}
Note that terms $A_{2,n}, A_{3,n,N}$ and $A_{5,n,N}$ are the same as in NPP,
$A_{4,n,N}$ is a minor improvement of the corresponding term in NPP, and
$A_{6,n,J}$ reflects the ``convergence rate'' in (\ref{Y_W}). Term $A_{1,n} $ of NPP (which does not appear
in our bounds) is ``absorbed'' in the term $\inf_{1\le  J\le n} A_{7,J}$ in the bounds (i)-(iii), below,
due to a somewhat a different approximation (see (\ref{ZN})).

\begin{prop}\label{berry}
Let the assumptions of Theorem \ref{tlcgauss} (iii) prevail, with $f_{k,n}\equiv f$ for all $1\leq k \leq n$, $n\in \N^*$ where $f:\R^\nu\to \R $ is a Lipschitz function % satisfying  $\E f^2\big ({\mbf X}\big)<\infty$
with
%and there exists $| \partial f|_\infty>0$ such that
$|f(x)-f(y)|\leq | \partial f |_\infty \, | x-y|$ for all $x\neq y \in \R^\nu$.
%Assume the additional condition to (\ref{Y_W}): for any $J \in \N^*$,
%there exists a sequence $(u_J(n))_{n\in \N}$ such as
%\begin{equation}
%\max_{-J\leq j \leq J} \Big \{\sup_{0\leq \tau \leq 1} \E \big|{\mbf X}_{\! n}([n\tau]+j)-{\mbf W}_{\! \tau}(j) \big|^2 \Big \}\leq u_J(n)\quad %\mbox{for any $n\in \N^*$ and}\quad u_J(n) \limiten 0.    \label{Y_W2}
%\end{equation}
Define  $S_n:=n^{-1/2} \sum_{t=1}^n f\left({\mbf X}_{\! n}(t)\right)$ and let $S$ be a zero-mean Gaussian random variable with a variance $\sigma_S^2:=\int_0^1  \sum_{j\in {\Z}} {\rm Cov}\big(
f({\mbf W}_\tau(0)), f({\mbf W}_\tau(j))\big)\, {\rm d}\tau <\infty$. Then:
\begin{enumerate}
\item[(i)] For any function $h$ twice continuously differentiable with bounded second derivative and for every $n>K$,
\begin{eqnarray}
\label{breuer1}
\hskip-.5cm \Big |\E \big [ h\big (S_n \big)\big ] -\E \big [ h\big ( S\big ) \big ] \Big | \leq |h''|_\infty \, \Big(\inf_{N\geq m}\big \{A_{2,N}+ A_{3,n,N}+ A_{4,n,N}+A_{5,n,N}\big \} + \inf_{1\le J\le n}\big \{ A_{6,n,J}+A_{7,J}\big \} \Big ). \qquad
\end{eqnarray}
\item[(ii)] For any Lipschitz function $h$, and for every $n>K$,
\begin{eqnarray}
\nonumber  && \hskip-1.5cm \Big |\E\big [h\big (S_n\big ) \big ]- \E \big [ h\big (S\big )\big ]\Big |\ \leq \
|h'|_\infty\, \Big  ( \frac 2 {\sigma_S} \, \inf_{1 \le J\le n}\Big \{ A_{6,n,J}+A_{7,J}\Big \}    \\
\label{mainestimate2} && \hspace{.2cm} + \, \inf _{N\geq m} \Big   \{ \Big (\frac 1 {2 \sigma_S}
+ \frac 1 {\big(  (2K+\nu^m)\E[f^2({\mbf X})] \big )^{1/2}} \Big )  A_{2,N}
+    \frac{A_{3,n,N} + A_{4,n,N}+A_{5,n,N}}{\big ( \sum_{\ell=m}^N \sigma_{\ell,n}^2 \big )^{1/2}} \Big \}  \Big  ). \quad \qquad
\end{eqnarray}
\item[(iii)] For any $z\in\R$, and for every $n>K$,
\begin{eqnarray}
\nonumber   && \hskip-0.7cm \big|\mathbb{P}(S_n\leq z) ] -\mathbb{P}(S\leq z)\big|\ \leq \ \frac 2 {\sigma_S} \,  \Big  ( \frac 2 {\sigma_S} \, \inf_{1\le J\le n}\Big \{ A_{6,n,J}+A_{7,J}\Big \}    \\
\label{mainestimate3} && \hspace{1.1cm} + \inf _{N\geq m} \Big   \{ \Big (\frac 1 {2 \sigma_S}
+ \frac 1 {\big(  (2K+\nu^m)\E[f^2({\mbf X})] \big )^{1/2}} \Big )  A_{2,N}
+    \frac{A_{3,n,N} + A_{4,n,N}+A_{5,n,N}}{\big ( \sum_{\ell=m}^N \sigma_{\ell,n}^2 \big )^{1/2}} \Big \}  \Big  )^{1/2}.\quad \qquad 
\end{eqnarray}
\end{enumerate}
\end{prop}

\noindent {\it Proof of Proposition \ref{breuer}.}
Let us introduce a similar notation to NPP. Consider the Hilbert space $\HH = \R^{n\nu} $ with elements
$u = (u_{t,l}, 1\le t  \le n, 1\le l \le \nu ) \in \HH $ and the scalar product
$\langle u_{t,j}, u_{t', j'} \rangle_{\HH} := \E X^{(j)}_n (t) X^{(j')}_n (t') =  r_n^{(j,j')}(t,t').  $ The $\ell-$fold tensor product and the symmetrized
tensor product of $\HH$ are denoted by
$\HH^{\otimes  \ell} $ and $\HH^{\odot  \ell}, $ respectively.
Let $\L^2({\mathfrak X}_n)$ denote the space of r.v.'s subordinated to the
Gaussian vector  ${\mathfrak X}_n := ({\mbf X}_{\! n}(t))_{1\le t \le n} $. Any element $\xi \in \L^2({\mathfrak X}_n) $ admits a
chaotic expansion $\xi = \sum_{\ell = 0}^\infty I_{(\ell)} (g_{(\ell)})$, where $g_{(\ell)} \in \HH^{\otimes  \ell}$ and
the linear mapping $I_{(\ell)}: \HH^{\otimes  \ell} \to \L^2({\mathfrak X}_n)$ satisfies
$I_{(\ell)}(g) = I_{(\ell)} ({\rm sym}(g)),$ $\E I^2_{(\ell)}(g) = \ell ! \|{\rm sym}(g)\|_{\HH^{\otimes  \ell}}$, and
$\E [I_{(\ell)}(g)I_{(\ell')}(g')] = 0, \, \ell \ne \ell',  \, g_{(\ell)} \in \HH^{\otimes  \ell}, \, g_{(\ell')} \in \HH^{\otimes  \ell'},$
where ${\rm sym}$ denotes the symmetrization operator.
In particular,  for any $t=1, \dots, n,\,  {\mbf k} \in \Z^\nu_+, |{\mbf k}| =: \ell $ we have $H_{\mbf k}({\mbf X}_{\! n}(t))
= I_{(\ell)} \big(g_{\ell}({\mbf k})\big) $, where
\begin{eqnarray*}\label{hermite1}
g_\ell({\mbf k})&:=&{\rm sym} \big(u^{\otimes k^{(1)}}_{t, 1}   \otimes \cdots \otimes
u_{t, \nu}^{\otimes k^{(\nu)}}\big) \
= \ \sum_{{\mbf v}\in \{1, \ldots, \nu\}^\ell } b({\mbf v}; {\mbf k}) \,u_{t,v_1} \otimes \cdots
\otimes u_{t,v_\ell}
\end{eqnarray*}
and where $b({\mbf v}; {\mbf k}) = {\rm sym} [\tilde  b({\mbf v}; {\mbf k})] $
%= (\ell !)^{-1} \sum_{ \pi \in \Pi_\nu} \tilde b(\pi {\mbf v}; {\mbf k})$
is the symmetrization of the function
$\{1, \dots, \nu\}^\ell \ni {\mbf v} = (v_1, \dots, v_\ell) \mapsto \tilde b({\mbf v}; {\mbf k}) :=  \prod_{r=1}^\nu
\1 \big(v_i = r, k_1 + \dots + k_{r-1} < i \le k_1 + \dots + k_r) $.
Thus, $S_n=n^{-1/2} \sum_{t=1}^n f\left({\mbf X}_{\! n}(t)\right) $ admits the chaotic expansion
$$
S_n=\sum_{\ell =m}^\infty I_{(\ell)}(g_\ell^n)\quad \mbox{with}\quad g_\ell^n := \frac{1}{\sqrt{n}}
\sum_{t=1}^n \sum_{{\mbf v}\in \{1, \ldots, \nu\}^\ell } b_\ell({\mbf v}) \,u_{t,v_1} \otimes \cdots
\otimes u_{t,v_\ell},
$$
where $ b_\ell ({\mbf v}) := \sum_{|{\mbf k}|= \ell} (J_f({\mbf k})/{\mbf k}!) b({\mbf v}; {\mbf k})$ depend only on
$f \in  \L^2({\mbf X}_{\! n}(t)) = \L^2({\mbf X})$ and satisfy
$\E f^2_{(\ell)}({\mbf X}) = \ell! \sum_{{\mbf v}\in \{1, \ldots, \nu\}^\ell } b^2_\ell({\mbf v}) $, as in NPP.
It is important that here the $g_\ell^n$'s are symmetric since the $b_\ell ({\mbf v})$'s are symmetric. Therefore
$\E I^2_{(\ell)}(g_\ell^n) = \ell!  \|g_\ell^n\|^2_{\HH^{\otimes \ell}}$.
Next,
for $N\geq m$ consider the truncated expansion
$$
S_{n,N}:= \sum_{\ell =m}^N I_{(\ell)}(g_\ell^n).
$$
Note that
\begin{eqnarray*}
\E S^2_{n,N}&=&\sum_{\ell = m}^N  \E I^2_{(\ell)}(g_\ell^n)  = \sum_{\ell = m}^N \ell! \|g_\ell^n\|^2_{\HH^{\otimes \ell}} \nonumber \\
&=&\frac{1}{n}\sum_{\ell = m}^N \ell!
\sum_{t,t'=1}^n \sum_{{\mbf v}, {\mbf v}'\in \{1, \ldots, \nu\}^\ell } b_\ell({\mbf v}) b_\ell({\mbf v}')\,
\langle u_{t,v_1} \otimes \cdots
\otimes u_{t,v_\ell}, u_{t',v'_1} \otimes \cdots
\otimes u_{t',v'_\ell} \rangle_{\HH^{\otimes \ell}} \nonumber \\
&=&\frac{1}{n}\sum_{\ell = m}^N \ell!
\sum_{t,t'=1}^n \sum_{{\mbf v}, {\mbf v}'\in \{1, \ldots, \nu\}^\ell } b_\ell({\mbf v}) b_\ell({\mbf v}')\,
\prod_{i=1}^\ell
r^{(v_i,v'_i)}_n (t,t'). \label{sigma!}
\end{eqnarray*}
Using $|r_n^{(j,j')}(t,t')|\leq \theta(t-t')$  similarly as in NPP we obtain
\begin{eqnarray}\label{SNn}
\Big|\E \big [h\big (S_n\big) \big]-\E \big [h\big (S_{n,N}\big)\big]\Big| \leq \frac 3 2 \, (2K + \nu^m \theta )\,  |h''|_{\infty} \,\Big ( \E[f^2({\mbf X})]\sum_{\ell=N+1}^\infty \E[f^2_{(\ell)}({\mbf X})]\Big )^{1/2}\leq \frac 3 4 \,|h''|_{\infty} \, A_{2,N}.
\end{eqnarray}
For $N\geq m$, let $Z_{n,N}$ be a
centered Gaussian random variable with variance $\E S^2_{n,N} =  \sum_{\ell=m}^N \sigma^2_{\ell,n}$, with
$\sigma^2_{\ell,n} $
% = \sum_{t=-n}^n  \Cov\big (f_{(\ell)}\left({\mbf
%X}_{\! n}(0) \right), f_{(\ell)}\left({\mbf
%X}_{\! n}(t) \right) \big)$
defined in (\ref{sigmam}). (Note that the last variance is slightly different from the variance of $Z_N$ in (NPP, sec. 4.2).)
Let $D$ denote the Malliavin derivative in $\L^2({\mathfrak X}_n)$, see NPP.
Using $\ell^{-1} \E \| D I_{(\ell)}(g_\ell^n)\|^2_{\HH} = \ell! \|g_\ell^n\|^2_{\HH^{\otimes \ell}}
= \sigma^2_{\ell, n}$, see (\ref{sigma!}), as in (NPP, (4.46)) we obtain
%The inequality (NPP, (4.45) becomes
\begin{eqnarray}
\nonumber \Big| \E \big [h\big ( Z_{n,N} \big) \big ]-\E \big [h\big ( S_{n,N}\big) \big ] \Big|
&\leq & \frac 1 2 \, |h''|_\infty \,\sum_{\ell, \ell'=m}^N \big \|\delta_{\ell \ell'} \sigma^2_{\ell,n} -
\ell^{-1} \langle D I_{(\ell)}(g_\ell^n), D I_{(\ell')}(g_{\ell'}^n)\rangle_{\HH}   \big \| _{\L^2(\P)} \nonumber \\
&\le&|h''|_\infty \big(A_{3,n,N} + A_{4,n,N} + A_{5,n,N}\big).
\label{ZN}
%&\leq & \|h''\|_\infty \, \big(A_{4,n,N}+A_{5,n,N}\big).
\end{eqnarray}
Next, using (NPP, (3.39))
\begin{eqnarray*}
\nonumber \Big| \E \big [h\big ( Z_{n,N} \big) \big ]-\E \big [h\big ( S\big) \big ] \Big|
&\leq &\frac 1 2 \,  |h''|_\infty \, \Big| \sum_{\ell=m}^N \sigma^2_{\ell,n} -\sigma_S^2\Big| \
\leq \ \frac 1 2 \,  |h''|_\infty \Big( \big|\sigma^2_{n} -\sigma_S^2\big| + \Big|\sigma^2_n -  \sum_{\ell=m}^N \sigma^2_{\ell,n} \Big|\Big).
%\label{ZS}
\end{eqnarray*} 
To estimate the difference $\sigma^2_{n} -\sigma_S^2$,
we use an interpolation identity from Houdr\'e {\it et al.} (1998). Let
$({\mbf X}_1, {\mbf X}_2), ({\mbf W}_1, {\mbf W}_2)$  be two $(2\nu)-$dimensional Gaussian vectors with zero means
and respective covariance matrices $\E [{\mbf X}_i {\mbf X}_i^\intercal]
= \E [{\mbf W}_i {\mbf W}_i^\intercal] = I, \, i=1,2, \, \E [{\mbf X}_1 {\mbf X}_2^\intercal] = \Sigma_1,\,
\E [{\mbf W}_1 {\mbf W}_2^\intercal] = \Sigma_0 $. For $\alpha \in [0,1]$ let
$({\mbf X}_{1\alpha}, {\mbf X}_{2\alpha})$ denote  the ``interpolated''  Gaussian vector with
zero mean and   $\E [{\mbf X}_{i\alpha} {\mbf X}_{i\alpha}^\intercal] = I, \, i=1,2, \,
 \E [{\mbf X}_{1\alpha} {\mbf X}_{2\alpha}^\intercal] = (1-\alpha) \Sigma_0 + \alpha \Sigma_1 $.
Let $ f \in \L^2({\mbf X})$ be a real  function satisfying the conditions of Proposition \ref{breuer}.
Then from (\cite{Houd}, (1.1), (1.3)) we obtain
\begin{eqnarray}
\big|\Cov\big (f\left({\mbf
X}_1 \right), f\left({\mbf
X}_2 \right)\big) - {\rm Cov}\big(
f({\mbf W}_1), f({\mbf W}_2)\big)\big|&=&\Big|
\int_0^1 \E \left[\partial f({\mbf X}_{1 \alpha})^\intercal  (\Sigma_1 - \Sigma_0) \partial f({\mbf X}_{2 \alpha}) \right]
\,{\rm d}\alpha \Big| \nonumber \\
&\le&|\partial f|_\infty^2 \, \|\Sigma_1 - \Sigma_0\|,  \label{Hin}
\end{eqnarray}
where $\partial f = (\partial f/\partial x^{(1)}, \dots, \partial f/\partial x^{(\nu)})^\intercal \in \R^\nu.  $
Let $F_n(\tau) := \sum_{t'=1}^n  \Cov\big (f\left({\mbf
X}_{\! n}([n\tau]) \right), f\left({\mbf
X}_{\! n}(t') \right)\big), \, \tau \in [0,1]    $  so that $ \sigma^2_{n} = \int_0^1 F_n(\tau) {\d} \tau. $
Using (\ref{Hin}), for $1\le J \le n$ we can write $\big |\sigma^2_{n}  -\sigma_S^2 \big| \le R_1(n,J) + R_2(n,J), $ where
\begin{eqnarray*}
R_1(n,J)
&:=&\int_0^1 \sum_{|j| \le J} \big|\Cov\big (f\big({\mbf
X}_{\! n}([n\tau]) \big), f\big({\mbf
X}_{\! n}([n\tau] +j) \big) - \Cov\big( f({\mbf W}_\tau(0)), f({\mbf W}_\tau(j))\big)\big| {\rm d} \tau \ \le \ 2 A_{6,n,J},\\
%&\le&|\partial f|^2_\infty \sup_{0\le \tau \le 1} \sum_{|j| \le J} \Big\| \E [{\mbf
%X}_{\! n}([n\tau]) {\mbf
%X}_{\! n}([n\tau] +j)^\intercal ] - \E [ {\mbf W}_\tau(0)) {\mbf W}_\tau(j)^\intercal ]\Big\|  =   A_{6,n,J}, \\
%&=:&|\partial f|^2_\infty u_{n,J}, \\
R_2(n,J)&\le&2\, \E \big [ f^2({\mbf X})\big ]   \,\nu^m \,  \sum_{|k|>J} \theta^m(k) \ = \ 2 A_{7, J}.
\end{eqnarray*}
We also have  %\begin{eqnarray*}
$\big |\sigma^2_{n}  -  \sum_{\ell = m}^N \sigma^2_{\ell,n} \big|
= \sum_{\ell = N+1}^\infty \sigma^2_{\ell,n} \ \le \  \frac 1 2 \,  A_{2,N}, $ as in (\ref{SNn}).
%\end{eqnarray*}
Therefore,  $\Big| \sum_{\ell=m}^N \sigma^2_{\ell,n} -\sigma_S^2\Big|
\le 2 A_{6,n,J}   + 2 A_{7,J} + \frac 1 2 \,  A_{2,N}, $ implying
\begin{eqnarray}
\Big| \E \big [h\big ( Z_{n,N} \big) \big ]-\E \big [h\big ( S\big) \big ] \Big|
&\leq &|h''|_\infty \Big(A_{6,n,J}   + A_{7,J} + \frac 1 4 \,  A_{2,N}  \Big) \qquad \text{for} \quad 1\le J \le n.
\label{ZS}
\end{eqnarray}
Finally combining (\ref{SNn}), ((\ref{ZN}), and (\ref{ZS}) results in
\begin{eqnarray*}
\Big |\E \big [ h\big (S_n \big)\big ] -\E \big [ h\big ( S\big ) \big ] \Big | &\leq &|h''|_\infty
\Big(A_{2,N} +  A_{3,n,N} + A_{4,n,N} + A_{5,n,N} + \inf_{1\le J \le n}(A_{6,n,J}   + A_{7,J})\Big) \\
&\leq &|h''|_\infty
\Big(\inf_{N \ge m} \big\{ A_{2,N} +  A_{3,n,N} + A_{4,n,N} + A_{5,n,N}\big\}   +
\inf_{1 \le J \le n} \big\{ A_{6,n,J} + A_{7,J}\big\}\Big),
\end{eqnarray*}
proving the bound in (\ref{breuer1}).   ~\\

\noindent (ii) Following (NPP, proof of Theorem 2.2-(2)) and the previous results, for a Lipschitz function $h$  we obtain:
\begin{eqnarray*}\label{SNn2}
\big | \E \big [h\big ( S_n\big) \big ]-\E \big [h\big ( S_{n,N}\big) \big ] \big |&\leq &|h'|_{\infty} \big(  (2K+\nu^m)\E[f^2({\mbf X})] \big )^{-1/2}  \, A_{2,N},\\
\big | \E \big [h\big ( Z_{n,N} \big) \big ]-\E \big [h\big ( S_{n,N}\big) \big ] \big |
\label{ZN2} &\leq &2 \,  |h'|_\infty \,\big ( \sum_{\ell=m}^N \sigma_{\ell,n}^2 \big )^{-1/2} \,  \big(A_{3,n,N}+  A_{4,n,N}+A_{5,n,N}\big)\\
\mbox{and}\quad \big | \E \big [h\big ( Z_{n,N} \big) \big ]-\E \big [h\big ( S\big) \big ] \big |& \leq &
\frac {|h'|_\infty}{\sigma_S}\,  \Big ( \frac 1 2 \,  A_{2,N}+ \inf_{1\le J \le n}(A_{6,n,J}+A_{7,J}) \Big )
\end{eqnarray*}
and therefore (\ref{mainestimate2}) is established. \\
~\\
(iii) Bound  (\ref{mainestimate3}) is obtained exactly as in  (NPP, proof of Theorem 2.2-(3)).  \hfill $\Box $ \\

\section{Applications of Lemma \ref{lemgauss} and Theorem \ref{tlcgauss}}\label{Applications}
\subsection{Application to the IR statistic}

This application was developed in Bardet and Surgailis (2011, 2012).
Let $(X_t)_{t\in [0,1]}$ be a  continuous time Gaussian process with zero mean and generally nonstationary increments  locally
resembling a fractional Brownian motion with Hurst parameter $H(t) \in (0,1)$.
%tangent process (which is a self-similar process with parameter $H(t)$) and
Consider the Increment Ratio (IR) statistic
\begin{eqnarray*}
R^{2,n} (X) &:=&\frac{1}{n-2} \sum_{k=0}^{n-3}
\frac{\big|\Delta^{2,n}_k X +  \Delta^{2,n}_{k+1} X \big|}
{|\Delta^{2,n}_k X| + |\Delta^{2,n}_{k+1} X| }, \label{R_p}
\end{eqnarray*}
with $\Delta^{2,n}_k X =X_{(k+2)/n}-2\, X_{(k+1)/n}+X_{k/n}$ and the convention $\frac{0}{0}:=1$. Let
$\sigma^2_{2,n}(k):=\E \Big[ \left(\Delta^{2,n}_k X\right)^2\Big] $ and
\begin{eqnarray*}
Y^{(1)}_{n}(k)&:=&\frac{\Delta^{2,n}_{k} X}{\sigma_{2,n}(k)}, \qquad Y^{(2)}_{n}(k) \ :=\ \frac{\Delta^{2,n}_{k+1} X}{\sigma_{2,n}(k)}.
\end{eqnarray*}
Then $R^{2,n} (X) = \frac 1 {n-2} \sum_{k=0}^{n-3} f\left({\mbf Y}_{\! n}(k)\right), \ f(x^{(1)}, x^{(2)}) :=
|x^{(1)}+x^{(2)}|/(|x^{(1)}|+|x^{(2)}|) $ can be written as the sum of nonlinear function $f $
of Gaussian vectors
${\mbf Y}_{\! n}(k) = (Y^{(1)}_n(k), Y^{(2)}_n(k)) \in
{\R}^2, \, 0 \le k \le n-3. $ These Gaussian  vectors can be standardized, leading to the expression
$R^{2,n} (X) = \frac 1 {n-2} \sum_{k=0}^{n-3} f_{n,k}\left({\mbf X}_{\! n}(k)\right) $ of the IR statistics  as the sum of some functions
$f_{n,k}$ of standardized Gaussian vectors  ${\mbf X}_{\! n}(k),  \, 0 \le k \le n-3. $ (It is easy to check that
the centered functions
$f_{n,k} - \E [f_{n,k}({\mbf X})]$ have the Hermite rank 2.)
If $(X_t)$ satisfies some additional conditions (specifying the decay rate of correlations of increments and the convergence
rate to the tangent process), Theorem \ref{tlcgauss} can be applied to establish 
that $\sqrt n \big (R^{2,n} (X) - \int_0^1 \Lambda(H(t))\, \d t\big)\limiteloin {\cal N} (0,\sigma^2) $ with an explicit function $\Lambda$ and a variance $\sigma^2.$ An application of Lemma \ref{lemgauss} to bound the 4th moment
$\E (R^{2,n} (X) - \E R^{2,n} (X))^4 $
provides a crucial step in the proof of the almost sure consistency of the IR statistic, {\it i.e.} $ R^{2,n} (X) \limitepsn \int_0^1 \Lambda(H(t))\, \d t$. See Bardet and Surgailis (2011) for details. Local versions of the IR statistic for point-wise estimation of $H(t)$ are
developed in Bardet and Surgailis (2012). The study of  the asymptotic properties of these estimators in the last paper  is also based on %an application of
Theorem \ref{tlcgauss} and  Lemma \ref{lemgauss}.
%Such results can be applied to fractional Brownian motions but as well to multifractional Brownian motions (without stationary properties). More %details can be seen in Bardet and Surgailis (2011).

\subsection{A central limit theorem for functions of locally stationary Gaussian processes}

Using an adaptation of Dahlhaus and Polonik (2006, 2009), we will say that $(X_{t,n})_{1\leq t \leq n,\,n \in \N^*}$ is a locally stationary Gaussian process if
\begin{equation} \label{local}
X_{t,n}:=\sum_{j\in \Z} a_{t,n}(j)\, \varepsilon_{t-j}, %  \ = \ \sum_{s\in \Z} a_{t,n}(t-s)\, \varepsilon_s
\qquad \mbox{for all $1\leq t \leq n,\, n\in \N^*$,}
\end{equation}
where $(\varepsilon_k)_{k\in \Z}$ is a sequence of independent standardized Gaussian variables and for $1\leq t\leq n$, $n\in N^*$ the sequences $(a_{t,n}(j))_{j\in \Z}$ are such that there exist $K\geq 0$ and $\alpha< 1/2$ satisfying for all $n\in \N^*$ and $j\in \Z$,
\begin{eqnarray}\label{condi1}
\max_{1\leq t\leq n} |a_{t,n}(j)|\leq \frac K {u_j},\qquad \mbox{ with \ $u_j:=\max (1,|j|^{\alpha-1})$ \ for $j \in \Z$}
\end{eqnarray} 
and such that there exist functions $\tau \in (0,1] \mapsto a(\tau,j)\in \R$ satisfying
the following conditions:
\begin{eqnarray}\label{conda}
\sup_{\tau \in (0,1]}|a(\tau,j)|&\leq&  \frac K {u_j}, \qquad \forall \, j \in \Z,\\
\label{cond3} \mbox{and}\qquad \sup_{\tau \in (0,1]}\, \max_{|[n \tau] - k| \le L} |(a_{k,n}(j) - a(\tau, j)\big|  \ &\to & \ 0,
\qquad \forall \, j \in \Z, \quad \forall \, L >0. 
\end{eqnarray}
For $\tau \in (0,1]$
introduce a stationary  Gaussian process
$$
W_{\tau}(t):= \sum_{j\in \Z} a(\tau,j)\, \varepsilon_{t-j}, \qquad t \in \Z.
$$
with spectral density $g_\tau(v) = |\hat a(\tau,v)|^2, \ \hat a(\tau, v) :=  (2\pi)^{-1/2} \sum_{j \in \Z} \e^{- \i j v} a(\tau, j), \,
v \in [-\pi, \pi]$. Let
$$
{\mbf Y}_{\!n}(k)\ :=\  \big (X_{k+1,n},\dots,X_{k+\nu,n})^\intercal, \qquad {\mbf W}_{\!\tau}(j)\ := \ \big (W_{\tau}(j+1),\dots, W_{\tau}(j+\nu) \big)^\intercal.
$$
Note $({\mbf W}_{\!\tau}(j))_{j \in \Z} $ is a $\R^\nu-$valued stationary Gaussian process. Let
$$
\Sigma_{k,n}\ := \  \E [ {\mbf Y}_{\! n}(k) {\mbf Y}_{\! n}(k)^\intercal], \qquad
\Sigma_\tau :=  \E [{\mbf W}_{\!\tau}(0) {\mbf W}_{\!\tau}(0)^\intercal ].
$$

\begin{prop} \label{localCLT} In addition to (\ref{local}) - (\ref{cond3}), assume that
\begin{equation} \label{bddS1}
\sup_{\tau \in (0,1]} \| \Sigma^{-1}_\tau \| < \infty.
\end{equation}
Let $f_{k,n} \in  \L^2_0({\mbf Y}_{\!n}(k)), \,  1 \le k \le  n,  n \ge 1$
% with ${\mbf Z}$ a standardized Gaussian vector $\R^d$-valued,
be a triangular array of functions all having a generalized Hermite rank at least $m > 1/(1-2\alpha) $. Let there exists a $ \L^2_0({\mbf X})-$valued
continuous function $\tilde \phi_\tau, \tau \in (0,1]$ such that relation (\ref{rank22}) holds, with
$\tilde f_{k,n}({\mbf x}) :=f_{k,n}( \Sigma^{1/2}_{k,n} {\mbf x}). $
Then the CLT of (\ref{CLTbis}) holds, with
\begin{eqnarray}\label{sigmaloc}
\sigma^2 :=\int_0^1 \d \tau \, \sum_{j\in \Z} \E \big [ \phi_\tau\big ({\mbf W}_{\!\tau}(0)\big)\, \phi_\tau\big({\mbf W}_{\!\tau}(j) \big)\big ]
\end{eqnarray}
and   $\phi_\tau ({\mbf x}) := \tilde \phi_\tau (\Sigma^{-1/2}_\tau {\mbf x})$
defined as in Corollary \ref{tlcgauss2}.

\end{prop}

\noindent {\it Proof.} We apply Corollary \ref{tlcgauss2}. Let us first check
\begin{equation} \label{approxS1}
\sup_{\tau \in (0,1]}\|\Sigma_{[n \tau],n} - \Sigma_\tau \| \limiten   0.
\end{equation}
We have
\begin{eqnarray*}
\big|\sigma_{[n\tau],n}(p,q) - \sigma_\tau (p,q)\big|&=&
\Big|\sum_{j \in \Z} \big(a_{[n \tau] +p,n}(p+j) a_{[n \tau]+ q,n} (q +j) - a(\tau, p+j) a(\tau, q +j)\big)\Big| \ \le \ T_{n,J} + T''_{n,J},
\end{eqnarray*}
where
\begin{eqnarray*}
T'_{n,J}&:=&2K^2\sum_{|j| > J}u_{p+j} u_{q+j}, \qquad T''_{n,J} :=  \sum_{|j| \le J} \big|a_{[n \tau] +p,n}(p+j) a_{[n \tau]+ q,n} (q +j) - a(\tau, p+j) a(\tau, q +j)\big)\big|
\end{eqnarray*}
according to (\ref{condi1}) and (\ref{conda}). Clearly, $T'_{n,J}$ can be made arbitrarily small by choosing $J$ large enough. Then
for any $J < \infty$ fixed, we have that $\sup_{\tau \in (0,1]} T''_{n,J} \to 0 $ according to assumption (\ref{cond3}). This proves
(\ref{approxS1}). In a similar way, one verify that for any $\tau \in (0,1], \, j, j' \in \Z$,
$\|\E [ {\mbf Y}_{\! n}([n\tau] +j) {\mbf Y}_{\! n}([n\tau] + j')^\intercal] -
\E [{\mbf W}_{\!\tau}(j) {\mbf W}_{\!\tau}(j')^\intercal ]\| \limiten   0 $ implying condition (\ref{Y_W}). The dominating
condition (\ref{rho_dom}) on cross-covariances is ensured by (\ref{condi1}) and the fact that $(1-2\alpha) m > 1$. The remaining
conditions of  Corollary \ref{tlcgauss2} are trivially satisfied. \hfill   $\Box $ \\

\begin{rem} {\rm
Dahlhaus and Polonik (2006, 2009) discussed  the short-memory case $(a_{t,n}(j))_{j\in \Z}\in \ell^1, \, 1\leq t \leq n$  only. On the
other hand,
%was considered and for any $1\leq t \leq n$, the sequence $(a_{t,n}(j))_{j\in \Z}\in \ell^1$.
condition (\ref{condi1}) allows for the long-memory case $(a_{t,n}(j))_{j\in \Z}\in \ell^2, \sum_{j\in \Z} |a_{t,n}(j)| = \infty$.
%$ is only required.
The last case is also discussed in Roueff and von Sachs (2010), where similar conditions as (\ref{condi1}) and (\ref{conda}) are provided in
spectral terms. It is not clear whether condition (\ref{cond3}) allows for jumps of the parameter curves $\tau \mapsto a (\tau, \cdot)$
as in Dahlhaus and Polonik (2006, 2009), in particular, for abrupt changes of the memory intensity of Gaussian process (\ref{local}). See also
Lavancier {\it et al.} (2011) for a related class of nonstationary moving average processes with long memory.
%spectral density.
%However the property of local stationarity is more general in Dahlhaus and Polonik (2006, 2009) because the parameter curves are allowed to have jumps %and Conditions (\ref{cond3}) and (\ref{cond4}) are replaced by
%$$
%\sup_{\scriptsize  \begin{array}{c} 0\leq x_0<\ldots<x_m\leq 1\\
%m\in \N^*\end{array}}  \sum_{k=1}^m |a(x_k,j)-a(x_{k-1},j)|\leq \frac K {u_j} \quad \mbox{and}  \quad \sup_{j\in \Z}\sum_{t=1}^n  |a_{t,n}(j)-a(\frac %t n, j)|\leq K.
%$$
}
\end{rem}

\begin{rem} {\rm Note that
${\mbf x}^\intercal \Sigma_\tau {\mbf x} = \int_{-\pi }^\pi g_\tau(v) \big| \sum_{j=1}^\nu \e^{\i j v} x^{(j)} \big|^2 \d v $
for any ${\mbf x} = (x^{(1)}, \dots, x^{(\nu)})^\intercal \in \R^\nu$. Therefore condition
$\inf_{v \in [-\pi, \pi], \tau \in (0,1]} g_\tau (v) \ge \gamma > 0$ on the spectral density of $(W_{\tau}(t))$
implies condition (\ref{bddS1}), since
${\mbf x}^\intercal \Sigma_\tau {\mbf x}
\ge c|{\mbf x}|^2, \, c := 2\pi \nu \gamma >0 $.

}
\end{rem}

\begin{rem} {\rm For stationary Gaussian long memory process,  condition  $m(1-2\alpha)>1$ was first obtained in Taqqu (1975).
Proposition \ref{localCLT} can be applied to prove the asymptotic normality of various
statistics of  locally stationary processes, see, e.g., Roueff and von Sachs (2010).
%including
%variance, covariance,..., of locally %stationary long memory processes. \\

}
\end{rem}

\medskip

\noindent {\bf Acknowledgment.} The authors are grateful to two
anonymous referees for valuable suggestions and comments
%their corrections and suggestions
that helped to  improve the original version
of the paper.


\begin{thebibliography}{99}
%\item{Adler, R.} (1981) {\em Geometry of Random Fields.} Wiley, New York.

\bibitem{Arcones, M.A.} Arcones, M.A. (1994) Limit theorems for nonlinear
functionals of a stationary Gaussian sequence of vectors.
{\em Ann. Probab.} {\bf 22}, 2242--2274.

\bibitem{bar2008} Bardet, J.-M., Doukhan, P., Lang G. and Ragache, N. (2008)
The standard Lindeberg method applied to weakly dependent processes. {\em ESAIM Probability
and Statistics} {\bf 12}, 154--172.

\bibitem{barsur} Bardet, J.-M. and Surgailis, D. (2011) Measuring the roughness of random paths by increment ratios. {\em Bernoulli}, {\bf 17}, 749--780.
%\bibitem{Ayache} Ayache, A. and L\'evy V\'ehel, J. (2004) On the identification of the pointwise H\"older exponent
%of the generalized multifractional Brownian motion. {\em Stochastic
%Process. Appl.} {\bf 111}, 119--156.


\bibitem{barsurmbm} Bardet, J.-M. and Surgailis, D. (2012) A new nonparametric estimator of the local Hurst function
of multifractional processes. Preprint.
%multifractal function of a Gaussian process. {\em Statist. Probab.
%Lett.} {\bf 39}, 337--345.

%\item{Benassi, A., Cohen, S. and Istas, J.} (2004) On roughness indices for fractional fields.
%{\em Bernoulli} {\bf 10}, 357--373.

\bibitem{Breuer} Breuer, P. and Major, P. (1983) Central limit theorems
for nonlinear functionals of Gaussian fields. {\em J.  Multivariate
Anal.} {\bf 13}, 425--441.

\bibitem{Chambers} Chambers, D. and Slud, E. (1989) Central limit theorems for
nonlinear functional of stationary Gaussian process. {\em Probab.
Th. Rel. Fields} {\bf 80}, 323--349.

\bibitem{Clem} Coulon-Prieur, C. and Doukhan, P. (2000)
A triangular central limit theorem under a new weak dependence condition. {\em Stat. Probab. Letters} {\bf 47} 61--68.


\bibitem{Csorgo} Cs\"org{\H o}, M. and Mielnichuk, J. (1996) The empirical process of a short-range dependent stationary
sequence under Gaussian subordination. {\em Probab. Th. Rel. Fields} {\bf 104},  15--25.


\bibitem{Dahl1} Dahlhaus, R. and Polonik, W. (2006) Nonparametric quasi-maximum likelihood estimation for Gaussian locally stationary processes.
{\em Ann. Statist.} {\bf 34}, 2790--2824.

\bibitem{Dahl2} Dahlhaus, R. and Polonik, W. (2009) Empirical spectral processes for locally stationary time series. {\em Bernoulli} {\bf 15},  1--39.

%\bibitem{Dave }
%Davidson, J. (1992)
%A central limit theorem for globally nonstationary near-epoch dependent functions of mixing
%processes.  {\em Econometric Theory} {\bf  8 }, 313--329.

\bibitem{Dede} Dedecker, J. and Merlev\` ede, F. (2002) Necessary and sufficient conditions for the conditional central limit theorem. {\em Ann. Probab.} {\bf 30}, 1044--1081

\bibitem{Dob} Dobrushin, R. L. and Major, P. (1979)
Non-central limit theorems for nonlinear functionals of Gaussian fields.
{\em Z. Wahrsch. Verw. Gebiete} {\bf 50}, 27--52.

\bibitem{Giraitis} Giraitis, L. and  Surgailis, D. (1985) CLT and other limit
theorems for functionals of Gaussian processes. {\em Z. Wahrsch.
verw. Gebiete} {\bf 70}, 191--212.

%\bibitem{Guyon} Guyon, X. and Le\`on, J. (1989) Convergence en loi
%des H-variations d'un processus gaussien stationnaire. {\em Ann.
%Inst. Poincar\'e} {\bf 25}, 265--282.

\bibitem{Guo} Guo, H. and Koul, H.L. (2008) Asymptotic inference in some
heteroscedastic regression models with long memory design and errors. {\em Ann. Statist.} {\bf 36}, 458--487


\bibitem{Guyon} Guyon, X. and Le\`on, J. (1989) Convergence en loi
des H-variations d'un processus gaussien stationnaire. {\em Ann.
Inst. Poincar\'e} {\bf 25}, 265--282.


\bibitem{Houd} Houdr\'e, C., P\'erez-Abreu, V. and Surgailis, D. (1998)
Interpolation, correlation identities, and inequalities for infinitely divisible variables.
{\em J. Fourier Anal. Appl.} {\bf 4}, 651--668.

\bibitem{Jac} Jacod, J. and Shiryaev, A.B. (1987) {\em Limit theorems for stochastic processes}. Springer-Verlag, Berlin.

%\bibitem{Jong }de Jong, R. M. (1997)
%Central limit theorems for dependent heterogeneous random variables.
%{\em Econometric Theory} {\bf 13}, 353--367.


\bibitem{koul}  Koul, H.L. and Surgailis, D. (2002)  Asymptotic expansion of the
empirical process of long memory moving averages, in: H. Dehling, Th. Mikosch and M. S{\o}rensen (Eds.),
{\it Empirical Process Techniques for Dependent Data}, Birkh\"auser, Boston, pp.213--239.



\bibitem{Lav} Lavancier, F., Leipus, R., Philippe, A. and Surgailis, D. (2011) Detection of non-constant long memory parameter.
Preprint.



\bibitem{npp} Nourdin, I., Peccati, G. and Podolskij, M. (2011) Quantitative Breuer-Major Theorems. {\em Stochastic Process. Appl. } {\bf 121}, 793--812.

\bibitem{Peligrad} Peligrad, M. and Utev, S. (1997)
Central limit theorem for linear processes.
{\em Ann. Probab.} {\bf 25} 443-456

\bibitem{Rio} Rio E. (1995) About the Lindeberg method for strongly mixing sequences. {\em  ESAIM Probability
and Statistics} {\bf 1}, 35--61.

\bibitem{rvs} Roueff, F. and von Sachs, R. (2010) Locally stationary long memory estimation. {\em Stochastic Process. Appl. } {\bf 121}, 813--844.

\bibitem{Sanchez} Sanchez de Naranjo, M.V. (1993)
Non-central limit theorems for nonlinear functionals of $k$ Gaussian fields. {\em
J. Multivariate Anal.} {\bf 44}, 227--255.

\bibitem{Sousou} Soulier, Ph. (2001)  Moment bounds and central limit theorem for
functions of Gaussian vectors. {\em Statist. Probab. Lett.} {\bf
54}, 193--203.

\bibitem{Surgailis} Surgailis, D. (2000) Long-range dependence and Appell rank.
{\em Ann. Probab.} {\bf 28}, 478--497.

\bibitem{Taqqu0} Taqqu, M.S.  (1975) Weak convergence to the fractional Brownian
motion and to the Rosenblatt process. {\em  Z. Wahrsch. verw. Gebiete} {\bf 31},  287--302.

\bibitem{Taqqu1} Taqqu, M.S. (1977)  Law of the iterated logarithm for sums of
non-linear functions of Gaussian variables that exhibit a long range
dependence. \ {\em Z. Wahrsch. verw. Gebiete} {\bf 40}, 203--238.

\bibitem{Taqqu2}  Taqqu, M.S. (1979)
Convergence of integrated processes of arbitrary Hermite rank.
{\em Z. Wahrsch. verw. Gebiete} {\bf 50}, 53--83.
\end{thebibliography}
\end{document}